\documentclass[oneside,11pt,reqno]{amsart}
\usepackage{subcaption}
\usepackage{amsmath, amscd}
\usepackage[utf8]{inputenc}
\usepackage{mathtools}
\usepackage{amssymb}
\usepackage{amsthm}
\usepackage{thmtools}
\usepackage{graphicx}
\usepackage{tikz-network}
\usepackage{mathrsfs}
\usepackage{stmaryrd}
\usepackage{enumitem}
\usepackage{amsfonts}
\usepackage{amssymb}
\usepackage{amsxtra}
\usepackage{tikz-cd}
\usepackage{hyperref}
\usepackage{tikz}  

\definecolor{mygreen}{HTML}{278444}

\DeclareMathOperator{\lcm}{lcm}

\newtheorem{theorem}{Theorem}[section]
\newtheorem{corollary}{Corollary}[theorem]
\newtheorem{lemma}[theorem]{Lemma}
\newtheorem{proposition}[theorem]{Proposition}
\newtheorem{remark}[theorem]{Remark}
\usetikzlibrary{arrows.meta,
                quotes}

\theoremstyle{definition}
\newtheorem{definition}{Definition}[section]

\usepackage{amsmath}
\usepackage{amssymb}
\usepackage{graphicx, overpic}
\usepackage{caption}
\usepackage{subcaption}
\usepackage{pinlabel}
\usepackage{floatrow}
\usepackage{soul}
\usepackage{tikz}
\usetikzlibrary{decorations.pathreplacing}
\usetikzlibrary{decorations.markings}
\usetikzlibrary{shapes.geometric}
\usepackage{color}
\usetikzlibrary{shapes.gates.logic.US,trees,positioning, arrows}
\usetikzlibrary{knots}
 \usepackage{float}
 
 \usepackage{eucal}
\usepackage[all]{xy}
\usepackage{multicol}

\usepackage{mathtools}

\newtheorem*{prop-MC}{Proposition 4.28}

\renewcommand{\bar}[1]{\overline{#1}}

\renewcommand{\emptyset}{\varnothing}

\renewcommand{\implies}{\Rightarrow}

\newcommand{\showcomments}{yes}

\usepackage{ifthen}
\newsavebox{\commentbox}
{\ifthenelse{\equal{\showcomments}{yes}}%
{\footnotemark
    \begin{lrbox}{\commentbox}
    \begin{minipage}[t]{1.25in}\raggedright\sffamily\tiny
    \footnotemark[\arabic{footnote}]}
{\begin{lrbox}{\commentbox}}}%
{\ifthenelse{\equal{\showcomments}{yes}}%
{\end{minipage}\end{lrbox}\marginpar{\usebox{\commentbox}}}
{\end{lrbox}}}

\begin{document}

\title[Torsion-Free Lattices in Baumslag-Solitar Complexes]{Torsion-Free Lattices in Baumslag-Solitar Complexes}

\begin{abstract}
 This paper classifies the pairs of nonzero integers $(m,n)$ for which the locally compact group of combinatorial automorphisms, Aut$(X_{m,n})$, contains incommensurable torsion-free lattices, where $X_{m,n}$ is the combinatorial model for Baumslag-Solitar group $BS(m,n)$. In particular, we show that Aut$(X_{m,n})$ contains abstractly incommensurable torsion-free lattices if and only if there exists a prime $p \leq \gcd(m,n)$ such that either $\frac{m}{\gcd(m,n)}$ or $\frac{n}{\gcd(m,n)}$ is divisible by $p$. In all these cases we construct infinitely many commensurability classes. Additionally, we show that when  Aut$(X_{m,n})$ does not contain incommensurable lattices, the cell complex $X_{m,n}$ satisfies Leighton’s property.
\end{abstract} 

\author{Maya Verma}
\address{University of Oklahoma, Norman, OK 73019-3103, USA}
\email{Maya.Verma-1@ou.edu}

\maketitle

\section{Introduction}

In this paper, we examine torsion-free uniform lattices in the combinatorial automorphism group Aut$(X_{m,n})$. Here $X_{m,n}$ is a combinatorial model for the Baumslag-Solitar group BS$(m,n)$. This study can be seen as an extension of \cite{forester2022incommensurable}, where the existence of incommensurable uniform lattices in Aut$(X_{m,n})$ is established for certain pairs of integers $(m,n)$.
The primary goal of this paper is to address the question posed in \cite{forester2022incommensurable}:  for which pairs $(m, n)$, does Aut$(X_{m,n})$ contain incommensurable lattices?

The following result is proved in \cite{forester2022incommensurable}.
\begin{theorem}\label{Forester's thm}\emph{(Forester)}
    The group \emph{Aut}$(X_{d,dn})$ contains uniform lattices that are not abstractly commensurable if one of the following holds:
    \begin{enumerate}
        \item $\gcd(d,n) \neq 1$.
        \item $n$ has a non-trivial divisor $p \neq n$ such that $p < d$, or
        \item $n < d$ and $d  \equiv 1 \pmod{n}$.
    \end{enumerate}
\end{theorem}

We will classify for which pairs of integers $(m,n)$ the locally compact group Aut$(X_{m,n})$ of combinatorial automorphisms contains incommensurable lattices. We will write any pair of positive integers in the form $(dm,dn)$ where, $d=\gcd(dm,dn)$ (equivalently, $\gcd(m,n)=1$). The main result of this paper is the following theorem:

\begin{theorem} \label{main theorem}
The locally compact group \emph{Aut}$(X_{dm,dn})$, for $m$ and $n$ coprime, contains abstractly incommensurable torsion-free uniform lattices if and only if there exists a prime $p \leq d$ such that either $p \mid m$ or $p \mid n$. Furthermore, \emph{Aut}$(X_{dm,dn})$ contains infinitely many commensurability classes of lattices when $p \mid m$ or $p \mid n$.

\end{theorem}

Theorem \ref{main theorem} shows that many of the complexes $X_{m,n}$ serve as combinatorial models for incommensurable groups and this behavior is only known in a few other cases: products of locally finite trees \cite{MR1839489}, \cite{MR2341837}, \cite{MR2694733}, and also some examples of Dergacheva and Klyachko \cite{MR4545211}.
In fact, we will show that when $m$ and $n$ have no divisors less than or equal to $d$ then any two lattices in Aut$(X_{dm,dn})$ are commensurable up to conjugacy.

We will say that a space satisfies the \emph{Leighton property} if any two compact spaces covered by $X$ admit a common finite sheeted covering up to isomorphism. This definition is motivated by Leighton's theorem, which is equivalent to the fact that trees satisfy the Leighton property.

\begin{theorem} \cite{leighton}\emph{(Leighton’s theorem)}  Let $G_1$ and $G_2$ be finite connected graphs with a common
cover. Then they have a common finite cover.
\end{theorem}
Woodhouse established the Leighton property for a family of CAT$(0)$ cube complexes exhibiting symmetry and homogeneity similar to regular graphs \cite{MR4647679}, and for trees with fins \cite{MR4243770}. Shepherd, Gardam, and Woodhouse proved the Leighton property for ``trees of objects" $X$ when Aut$(X)$ has finite edge isometry groups \cite{MR4506536}. Additional variations on Leighton's theorem can be found in \cite{MR4464467}, \cite{MR4506536}, \cite{MR4647679}, \cite{MR4243770}, \cite{MR2727669}, \cite{MR4748179}, \cite{MR4545211}, \cite{MR1839489}, \cite{MR2694733}, \cite{MR2341837}. In the complementary case to Theorem \ref{main theorem}, when incommensurable lattices do not exist, we show that Leighton's property holds:

\begin{theorem} \label{Leighton property for X_{m,n}}
    The Baumslag-Solitar complex $X_{dm,dn}$ with $m$ and $n$ coprime has the Leighton property if and only if $m$ and $n$ have no divisors less than or equal to $d$.
\end{theorem}

Lastly, in addition to the result above, we have some new results about the invariant defined by Casals-Ruiz, Kazachkov, and Zakharov \cite{MR4359918}. In their work, the authors of solved the isomorphism problem for the following class of GBS groups by giving an isomorphism invariant vector (we will call it the CRKZ invariant) well-defined up to cyclic permutation;
 
$$BS(1,n^l) \vee  BS(n^{a_1}, n^{a_1}) \vee  BS(n^{a_2}, n^{a_2}) \vee \cdots \vee BS(n^{a_{k-1}}, n^{a_{k-1}})$$
for every $l \geq 1$, $n \geq 2$, $k \geq 2$, $0\leq a_1, a_2, \cdots, a_{k-1} \leq l-1$. This class of groups is denoted by $\mathcal{C}_{n,l}$. We show that this vector is a commensurability invariant up to scalar multiplication.

\begin{theorem}
    Suppose $G_1$ and $G_2$ are two groups in $\mathcal{C}_{n,l}$. Then $G_1$ is commensurable to $G_2$ if and only if $c_1\Vec{X}^l(G_1) = c_2\Vec{X}^l(G_2) $ for some $c_1, c_2 \in \mathbb{N}$. 
\end{theorem}
The complete statement of this result is given in Theorem \ref{5, 15 thm}.
\subsection{Methods} 
$X$ satisfies Leighton's property if and only if any two torsion-free lattices in Aut$(X)$ are commensurable up to conjugacy.

 Also, if there exists a prime $p \leq d$ such that either $p \mid m$ or $p \mid n$ then we will provide infinitely many examples of abstractly incommensurable lattices in Aut$(X_{dm,dn})$.  Therefore, to prove theorems \ref{main theorem}, and \ref{Leighton property for X_{m,n}}, it suffices to prove the ``if'' direction in both theorems.

Section \ref{construction of admissible cover for p-uni} is devoted to a general construction, which produces a prime power index subgroup of a $p$-unimodular GBS group. 

The converse in Theorem \ref{Leighton property for X_{m,n}} is proved in Section \ref{commensurability using leighton's property}. For $(dm,dn) \neq (1,1)$, if $X$ is a compact cell complex covered by $X_{dm,dn}$ then the fundamental group of $X$ is virtually a $p$-unimodular GBS group for all primes $1 \leq p \leq d$. We will use the result of Section \ref{construction of admissible cover for p-uni} and Leighton's graph theorem to construct isomorphic finite index subgroups of the fundamental groups of any two compact cell complexes covered by $X_{dm,dn}$. The proof for $X_{1,1}$ follows from the fact that Aut$(X_{1,1})$ is a discrete group and any two lattices in a discrete group are commensurable.

While proving the converse in Theorem \ref{main theorem}, without loss of generality we can assume that $p \mid m$. We split the proof into four cases, and provide examples of incommensurable lattices in Aut$(X_{dm,dn})$ for each case using different methods.

\begin{description}
     \item[\rm Case I] $p \mid d$;
 
 \item [\rm Case II] $p \nmid d$, $n =1$, and $m\neq p$; .
  \item[\rm Case III] $p \nmid d$, $n=1$, and $m=p$; 
    \item[\rm Case IV] $p \nmid d$ and $n > 1$.
\end{description}

Case I and II are proved in Section \ref{eg using DP}. In these cases, we employ depth profiles to construct incommensurable lattices in Aut$(X_{dm,dn})$. The depth profile is a commensurability invariant, developed by Forester in \cite{forester2022incommensurable}, taking the form of a subset of the natural numbers, depending on a choice of elliptic subgroup.

 In Section \ref{eg using crkz} we generalize the CRKZ invariant vector for a larger class of GBS groups denoted by $\mathcal{C}_{n,l}$, whose image under the modular homomorphism is generated by $1/n^l$ for some $l \geq 1$, and the index of an edge group in a vertex group is $n^i$. We will prove that the CRKZ invariant provides a commensurability invariant for certain GBS groups in  $\mathcal{C}_{n,l}$. For case III we will construct incommensurable lattices belonging in $\mathcal{C}_{n,l}$. 

Case IV is proved in Section \ref{eg using plateau}. The argument in this case depends on the notion of a $p$--plateau and slide equivalence of certain GBS groups.

\section{acknowledgement}
I would like to express my deepest gratitude to my advisor, Dr. Max Forester, for his unwavering support and invaluable guidance throughout every stage of this research paper. His insight and expertise have been instrumental in the development and completion of this paper. I am also grateful for the financial support provided by the National Science Foundation under grant number $1651963$ during the summer of 2023.

\section{preliminaries}
We refer to \cite{forester2022incommensurable} for more details on this section. 

A \textit{graph} $A$ consists of two sets $V(A)$ and $E(A)$, called the \textit{vertices} and \textit{edges} of $A$, respectively. It also includes an involution on $E(A)$, which send $e \in E(A)$ to $\overline{e} \in E(A)$, where $e \neq \overline{e}$, and maps $\partial_0, \partial_1: E(A) \shortrightarrow V(A)$, satisfying $\partial_1(e)=\partial_0(\overline{e})$.
For an edge $e$, the vertices $\partial_0(e)$ and $\partial_1(e)$ are called \textit{initial} and \textit{terminal} vertices of $e$, respectively. We say $e$ joins the initial vertex $\partial_0(e)$ to the terminal vertex $\partial_1(e)$.
An edge $e$ is called a $loop$ if $\partial_0(e)=\partial_0(\overline{e})$. For each vertex $v \in V(A)$, define $E_0(v)= \{e \in E(A): \partial_0(e)=v\}$.

A \textit{directed graph} is a graph $A$ together with a partition $E(A) = E^+(A) \sqcup E^{-}(A)$
that separates every pair $\{e, \overline{e}\}$. The edges in $E^+(A)$ are called \textit{directed edges}. For each $v \in V(A)$, we define 
$E^+_0(v)= \{ e \in E^{+}(A): \partial_0(e)=v\}$ and $E^-_0(v)= \{ e \in E^{-}(A): \partial_0(e)=v\}$.

A \emph{labeled graph} $(A, \lambda)$ is a finite graph with a label function $\lambda: E(A) \shortrightarrow (\mathbb{Z}-\{0\})$, hence each $e\in E(A)$ has a label $\lambda(e)$, which is a nonzero integer.

For a CW complex $X$, the \textit{topological automorphism group} (denoted as Aut$_{top}(X)$) consists of homeomorphisms of $X$ that preserve the cell complex structure of $X$. The \textit{combinatorial automorphism group}, denoted Aut$(X)$, is obtained by quotienting Aut$_{top}(X)$, where two automorphisms are considered the same if they induce the same permutation on the set of cells of $X$. For a connected and locally finite CW complex $X$, the combinatorial automorphism group Aut$(X)$ is locally compact.

In a locally compact group $G$, a discrete subgroup $H < G$ is called a \textit{lattice} if $G / H$ carries a finite positive $G$-invariant measure, and a \textit{uniform lattice} if $G / H$ is compact. A subgroup
$H$ in $G=$Aut$(X)$ is discrete if and only if every cell stabilizer $H_{\sigma}=\{ h \in H : h\sigma=\sigma\}$ is finite.  In this case, define the
covolume of $H$ to be:
$$\text{Vol}(X/H)=\sum_{[\sigma] \in cell(X/H)} 1/ |H_{\sigma}|  $$
where the sum is taken over a set of representatives of the $H$–orbits of cells of $X$. The next proposition follows from \cite[1.5-1.6]{basslubotzky},.

\begin{proposition}\label{freely and compactly}
    Let $X$ be a connected locally finite CW complex. Suppose that $G =\textup{Aut}(X)$ acts cocompactly on $X$ and let $H < G$ be a discrete subgroup. Then
    \begin{enumerate}
        \item $H$ is a lattice if and only if Vol$(X/H) < \infty$
        \item $H$ is a uniform lattice if and only if $X/H$ is compact.
    \end{enumerate}  
\end{proposition}
Note that, if $H$ is a torsion-free lattice in  Aut$(X)$ then every cell stabilizer is $H_{\sigma}$ is the trivial subgroup. Therefore  Vol$(X/H) < \infty$ iff $X/H$ is compact. Hence, the above Proposition implies that a torsion-free lattice is uniform.

\subsection{GBS group}  A \textit{generalized Baumslag--Solitar group} or a \textit{GBS group} is the fundamental group of a graph of groups where all edge and vertex groups are $\mathbb{Z}$. Any GBS group can be represented by a labeled graph $(A, \lambda)$, where the inclusion from the edge group $G_e$ to the vertex group $G_{\partial_0(e)}$ is given by multiplication by nonzero integer $\lambda(e)$.

The Baumslag-Solitar group, $BS(m,n)$ is the GBS group represented by the labeled graph with one vertex and one edge with labels $m$ and $n$.

For a labeled graph $(A, \lambda)$, a \emph{fiberd $2$-complex} denoted by $Z_{(A, \lambda)}$ is the total space of a graph of spaces in which vertex and edge spaces are circles. If $C_v$ is the oriented circle for $v \in V(A)$ and  $C_e$ is the oriented circle for $e \in E(A)$, and  $M_e$  is the mapping cylinder for the covering map $C_e \shortrightarrow C_{\partial_0(e)}$ of degree $\lambda(e)$, then 
$$Z_{(A, \lambda)}=\sqcup_{e\in E(A)} M_e \Big{/} \sim$$ 
where $M_e$ and $M_{\overline{e}}$ are identified along $C_e$ and $C_{\overline{e}}$ for each $e \in E(A)$, and all copies of $C_v$ are identified by identity map for each $v \in V(A)$.
The fundamental group of $Z_{(A, \lambda)}$ is the GBS group represented by $(A, \lambda)$.

\textbf{Notation:} We denote the GBS group defined by the labeled graph having one vertex and $k$ loops each with labels $m_i$ and $n_i$ by $\bigvee_{i=1}^{k}BS(m_i,n_i)$.

A labeled graph is called \textit{reduced} if every edge $e$ with $\lambda(e)= \pm1$  is a loop.

A \textit{G–tree} is a simplicial tree $X$ on which $G$ acts without inversions, i.e., if $g\in G$ fixes an edge, then it fixes every point on this edge. For a given  $G$-tree  $X$, an element $g \in G$ is called \textit{elliptic}  if it fixes a vertex, and \textit{hyperbolic} otherwise. Every hyperbolic element has a $g$-invariant line on which it acts via non-trivial translation. A subgroup $H< G$ is called elliptic if there exists a vertex  $v\in V(X)$ such that $hv=v$ for all $h\in H$.

A GBS group is called an \textit{elementary} GBS group if it is isomorphic to $\mathbb{Z}$, $\mathbb{Z} \times \mathbb{Z}$, the Klein bottle group or the union of infinitely ascending chain of infinite cyclic groups; otherwise, it is called a \textit{non-elementary} GBS group. If $G$ is a non-elementary GBS group, then 
any two $G$-trees produce the same set of elliptic (and hyperbolic) elements in $G$. 
Therefore, for a non-elementary GBS group, we can define the notion of elliptic (and hyperbolic) elements independently of its $G$-trees.

Let $G$ be a GBS group with a $G$-tree $X$ and a quotient labeled graph $(A, \lambda)$. Fix an elliptic element $a \in G$. Then, for any $g\in G$, we can find nonzero integers $m$ and $n$ such that $g^{-1}a^mg=a^n$.
The \textit{modular homomorphism} is a map $q: G \to \mathbb{Q}^*$ defined by  $q(g)=\frac{m}{n}$. This map is independent of the choice of elliptic element $a$. 
The restriction of the modular homomorphism to elliptic elements is a trivial map and $Q^*$ is abelian; therefore, it factors through $H_1(A)$. If $g\in G$ maps to $\alpha \in H_1(A)$, which is represented by a $1-$cycle $(e_1, e_2,\cdots, e_n)$, then
\begin{equation}\label{eq0}
    q(g) = \prod_{i=1}^n \frac{\lambda (e_i)}{\lambda (\overline{e}_i)}
\end{equation}
If $V$ is any non-trivial elliptic subgroup of $G$, then we have a formula:
\begin{equation}\label{q in terms of V}
|q(g)|=\frac{[V:V \cap V^g]}{[V^g: V \cap V^g]}
\end{equation}
 
\subsection{The $\mathbf{2-}$complex $\mathbf{X_{m,n}}$}

For $m,n >0$, we define $Z_{m,n}$ as the presentation $2-$complex corresponding to the group presentation $\langle a,t: t^{-1}a^mt=a^n \rangle$ for the Baumslag-Solitar group $BS(m,n)$. The complex  $Z_{m,n}$ consists of a single vertex, two edges labeled as $a$ and $t$, and a $2-$cell attached along the boundary word $t^{-1}a^mta^{-n}$. 
 We structure $Z_{m,n}$ as a cell complex by subdividing the $2$-cell into $\gcd(m,n)$ $2$-cells (see Figure \ref{Xmn} for an example). We will denote the universal cover of $Z_{m,n}$ by $X_{m,n}$, and the cell complex structure is inherited from $Z_{m,n}$. For any nonzero integers $m$ and $n$ define  $X_{m,n}$ to be  $X_{|m|,|n|}$.

\begin{figure}[!ht]
\begin{tikzpicture}[scale=0.9]
\small

\filldraw[fill=gray!15,thick] (1,1) rectangle (4,2);
\filldraw[fill=gray!15,thick] (5,1) rectangle (8,2);
\filldraw[fill=gray!15,thick] (11,1) rectangle (14,2);

\filldraw[fill=black,thick] (1,1) circle (.5mm);
\filldraw[fill=black,thick] (2.5,1) circle (.5mm);
\filldraw[fill=black,thick] (4,1) circle (.5mm);
\filldraw[fill=black,thick] (5,1) circle (.5mm);
\filldraw[fill=black,thick] (6.5,1) circle (.5mm);
\filldraw[fill=black,thick] (8,1) circle (.5mm);
\filldraw[fill=black,thick] (11,1) circle (.5mm);
\filldraw[fill=black,thick] (12.5,1) circle (.5mm);
\filldraw[fill=black,thick] (14,1) circle (.5mm);

\filldraw[fill=black,thick] (1,2) circle (.5mm);
\filldraw[fill=black,thick] (2,2) circle (.5mm);
\filldraw[fill=black,thick] (3,2) circle (.5mm);
\filldraw[fill=black,thick] (4,2) circle (.5mm);
\filldraw[fill=black,thick] (5,2) circle (.5mm);
\filldraw[fill=black,thick] (6,2) circle (.5mm);
\filldraw[fill=black,thick] (7,2) circle (.5mm);
\filldraw[fill=black,thick] (8,2) circle (.5mm);
\filldraw[fill=black,thick] (11,2) circle (.5mm);
\filldraw[fill=black,thick] (12,2) circle (.5mm);
\filldraw[fill=black,thick] (13,2) circle (.5mm);
\filldraw[fill=black,thick] (14,2) circle (.5mm);

\draw[thick,->] (1,1.54) -- (1,1.55);
\draw[thick,->] (4,1.54) -- (4,1.55);
\draw[thick,->] (5,1.54) -- (5,1.55);
\draw[thick,->] (8,1.54) -- (8,1.55);
\draw[thick,->] (11,1.54) -- (11,1.55);
\draw[thick,->] (14,1.54) -- (14,1.55);

\draw[thick,->] (1.84,1) -- (1.85,1);
\draw[thick,->] (3.34,1) -- (3.35,1);
\draw[thick,->] (5.84,1) -- (5.85,1);
\draw[thick,->] (7.34,1) -- (7.35,1);
\draw[thick,->] (11.84,1) -- (11.85,1);
\draw[thick,->] (13.34,1) -- (13.35,1);

\draw[thick,->] (1.5,2) -- (1.6,2);
\draw[thick,->] (2.5,2) -- (2.6,2);
\draw[thick,->] (3.5,2) -- (3.6,2);
\draw[thick,->] (5.5,2) -- (5.6,2);
\draw[thick,->] (6.5,2) -- (6.6,2);
\draw[thick,->] (7.5,2) -- (7.6,2);
\draw[thick,->] (11.5,2) -- (11.6,2);
\draw[thick,->] (12.5,2) -- (12.6,2);
\draw[thick,->] (13.5,2) -- (13.6,2);

\draw (9.52,1.46) node {$\dotsm$};

\draw (1.75,0.9) node[anchor=north] {$a$};
\draw (3.25,0.9) node[anchor=north] {$a$};
\draw (5.75,0.9) node[anchor=north] {$a$};
\draw (7.25,0.9) node[anchor=north] {$a$};
\draw (11.75,0.9) node[anchor=north] {$a$};
\draw (13.25,0.9) node[anchor=north] {$a$};
\draw (1.5,2.1) node[anchor=south] {$a$};
\draw (2.5,2.1) node[anchor=south] {$a$};
\draw (3.5,2.1) node[anchor=south] {$a$};
\draw (5.5,2.1) node[anchor=south] {$a$};
\draw (6.5,2.1) node[anchor=south] {$a$};
\draw (7.5,2.1) node[anchor=south] {$a$};
\draw (11.5,2.1) node[anchor=south] {$a$};
\draw (12.5,2.1) node[anchor=south] {$a$};
\draw (13.5,2.1) node[anchor=south] {$a$};
\draw (0.9,1.5) node[anchor=east] {$t$};
\draw (4,1.5) node[anchor=east] {$t_1$};
\draw (5.05,1.5) node[anchor=west] {$t_1$};
\draw (8.1,1.5) node[anchor=west] {$t_2$};
\draw (11,1.5) node[anchor=east] {$t_{k-1}$};
\draw (14.1,1.5) node[anchor=west] {$t$};

\end{tikzpicture}
\caption{The cell structure for $Z_{m,n}$ when $m=3k$,
  $n=2k$. 
  }\label{Zcells}
  \label{Xmn}
\end{figure}

The following theorem gives the sufficient condition for a GBS group to be a lattice in Aut$(X_{dm,dn})$. We will use this
theorem extensively throughout all sections to construct uniform lattices in Aut$(X_{dm,dn})$.

\begin{theorem}\cite{forester2022incommensurable}\label{labled graphs indeed define lattices}
    Let $G$ be the GBS group defined by labeled graph $(A, \lambda)$, and  suppose there is a directed graph structure $E(A)=E^+
    (A) \sqcup E^{-}(A)$ on $A$ such that
    \begin{enumerate}
        \item for every $v \in V(A)$, 
    $$\sum_{e \in E^+_0(v)} |\lambda(e)| =dm \: \; \text{and}
        \sum_{e \in E^-_0(v)} |\lambda(e)|=dn$$

        \item for every $e \in E^+(A)$, let $n_e=|\lambda(e)|$, $m_e=|\lambda(\overline{e})|$, and $k_e=\gcd(m_e,n_e)$; then
        $$n_e/k_e=n \text{ and } m_e/k_e=m.$$
    \end{enumerate}
    Then $G$ is  a lattice in Aut$(X_{dm,dn})$.
\end{theorem}

 The following theorem from \cite{forester2022incommensurable} provides the general description of a torsion-free uniform lattice within the combinatorial automorphism group Aut$(X_{dm,dn})$ as a GBS group, for $(m,n)\neq (1,1)$.
\begin{theorem}\label{prop4.6}
    \cite{forester2022incommensurable}  Suppose $m \neq n$ and let $G$ be a torsion-free group. Then $G$ is isomorphic to a uniform lattice in Aut$(X_{dm,dn})$ if and only if there exists a compact GBS structure $(A, \lambda)$ for $G$, a directed graph structure $E(A)=E^+(A) \sqcup E^-(A)$, and a length function $l: V(A) \sqcup E(A) \to \mathbb{N}$ satisfying $l(e)=l(\overline{e})$ for all $e \in E(A)$ such that the following holds. 
    \begin{enumerate}[label={(\arabic*)}]
    \item For every $v\in V(A)$: 

    \begin{align}
        \sum_{e \in E^+_0(v)} |\lambda(e)| =dm \: \; \text{and}
        \sum_{e \in E^-_0(v)} |\lambda(e)|=dn
    \end{align}

\item For every $e\in E^+(A)$,
\begin{align*}
    l(\partial_0(e) |\lambda(e)|=m l(e) \\
    l(\partial_1(e) |\lambda(\overline{e})|=n l(e)
\end{align*}

\item For every $v\in V(A)$, let $k_0(v)=\gcd(l(v), m)$ and $k_1(v)=\gcd(l(v), n)$; then there exist partitions

\begin{align*}
     E^+_0(v)=E^+_1 \sqcup \cdots \sqcup E^+_{k_0(v)} \: \; \text{and} \: \;
        E^-_0(v)=E^-_1 \sqcup \cdots \sqcup E^-_{k_0(v)}
\end{align*}
such that the sums $\sum_{e\in E^+_i} l(e)$ are all equal for all $i$, and the sums $\sum_{e\in E^-_j} l(e)$ are all equal for all $j$.
 
 \end{enumerate}
\end{theorem}

\begin{remark}\label{A is strongly connected}
By the proof of Proposition 4.6 in  \cite{forester2022incommensurable}, the conditions (1) and (2) in Theorem \ref{prop4.6} together imply that the directed graph is strongly connected. In particular, every directed edge is contained in a directed circuit.
    \end{remark}

\subsection{Deformation moves}Any two GBS trees for a non-elementary GBS group $G$ are related by an elementary deformation. That is, they are related by a finite sequence
of elementary moves, called elementary collapses and expansions.  There are also slide and induction moves, which can be expressed as an expansion followed by a collapse.

Collapse and expansion moves are as follows:
\newenvironment{pict}[2]%
	{\setlength{\unitlength}{1mm}
	\begin{center}
	\begin{picture}(#1,#2)
	\scriptsize}%
	{\end{picture}
	\end{center}

	\noindent}

\newcommand{\zindex}[3]{\put(#1,#2){\makebox(0,0){${#3}$}}}

\begin{pict}{90}{10}
\thicklines
\put(10,5){\circle*{1}}
\put(25,5){\circle*{1}}

\put(10,5){\line(1,0){15}}

\thinlines
\put(45,6.5){\vector(1,0){15}}
\put(60,2.2){\vector(-1,0){15}}
\zindex{52.5}{8.2}{\mbox{collapse}}
\zindex{52.5}{4}{\mbox{expansion}}

\put(25,5){\line(3,5){3}}
\put(25,5){\line(3,-5){3}}

\put(10,5){\line(-5,3){5}}
\put(10,5){\line(-5,-3){5}}

\put(80,5){\circle*{1}}
\put(80,5){\line(-5,3){5}}
\put(80,5){\line(-5,-3){5}}

\put(80,5){\line(3,5){3}}
\put(80,5){\line(3,-5){3}}

\scriptsize
\zindex{8.5}{8}{a}
\zindex{8.5}{2}{b}
\zindex{12}{6.5}{n}
\zindex{23}{6.5}{1}
\zindex{28.5}{7.5}{c}
\zindex{28.5}{2.5}{d}

\zindex{78.5}{8}{a}
\zindex{78.5}{2}{b}
\zindex{84.5}{7.5}{nc}
\zindex{84.5}{2.5}{nd}
\end{pict}

The induction move is as follows:

\begin{pict}{80}{10}
\thicklines
\put(6,5){\circle{10}}
\put(11,5){\circle*{1}}

\put(68,5){\circle{10}}
\put(73,5){\circle*{1}}

\thinlines
\put(11,5){\line(1,1){4}}
\put(11,5){\line(1,-1){4}}

\put(73,5){\line(1,1){4}}
\put(73,5){\line(1,-1){4}}

\scriptsize
\zindex{15}{7}{a}
\zindex{15}{3.2}{b}
\zindex{9}{6.7}{1}
\zindex{8}{3.3}{lm}

\put(31.5,5){\vector(1,0){17}}
\put(48.5,5){\vector(-1,0){17}}
\zindex{40}{7}{\mbox{induction}}

\zindex{78}{7}{la}
\zindex{78}{3}{lb}
\zindex{71}{6.7}{1}
\zindex{70}{3.3}{lm}
\end{pict}

There are two slide moves:
\begin{pict}{100}{11}
\thicklines
\put(75,3){\circle*{1}}
\put(90,3){\circle*{1}}
\put(75,3){\line(1,0){15}}
\put(90,3){\line(-1,2){4}}

\thinlines
\put(47.5,3){\vector(1,0){10}}
\zindex{52.5}{5}{\mbox{slide}}

\put(90,3){\line(5,3){5}}
\put(90,3){\line(5,-3){5}}
\put(75,3){\line(-5,3){5}}
\put(69,3){\line(1,0){6}}
\put(75,3){\line(-5,-3){5}}

\scriptsize
\zindex{77}{1.5}{m}
\zindex{88}{1.5}{n}
\zindex{86.5}{6}{ln}

\thicklines
\put(10,3){\circle*{1}}
\put(25,3){\circle*{1}}
\put(10,3){\line(1,0){15}}
\put(10,3){\line(1,2){4}}

\thinlines
\put(25,3){\line(5,3){5}}
\put(25,3){\line(5,-3){5}}
\put(10,3){\line(-5,3){5}}
\put(4,3){\line(1,0){6}}
\put(10,3){\line(-5,-3){5}}

\scriptsize
\zindex{12}{1.5}{m}
\zindex{23}{1.5}{n}
\zindex{9.5}{7}{lm}
\end{pict}

\begin{pict}{100}{10}
\thicklines

\put(82.5,5){\circle*{1}}
\put(87.5,5){\circle{10}}
\put(72.5,5){\line(1,0){10}}

\thinlines
\put(47.5,5){\vector(1,0){10}}
\zindex{52.5}{7}{\mbox{slide}}

\put(82.5,5){\line(-5,-3){4.5}}
\put(82.5,5){\line(-1,-4){1.2}}
\put(82.5,5){\line(-1,4){1.2}}

\scriptsize
\zindex{85.1}{3.5}{m}
\zindex{84.7}{6.5}{n}
\zindex{79.9}{6.8}{ln}

\thicklines

\put(17.5,5){\circle*{1}}
\put(22.5,5){\circle{10}}
\put(7.5,5){\line(1,0){10}}

\thinlines
\put(17.5,5){\line(-5,-3){4.5}}
\put(17.5,5){\line(-1,-4){1.2}}
\put(17.5,5){\line(-1,4){1.2}}

\scriptsize
\zindex{20.1}{3.5}{m}
\zindex{19.7}{6.5}{n}
\zindex{14.5}{6.8}{lm}
\end{pict}

We will use the following theorem by Forester \cite{MR2276234} in section \ref{eg using plateau}, which provides a sufficient condition for two reduced labeled graphs to be related only through slide moves. 

\begin{proposition}\label{reduced graph related by slide move}
    Suppose $(A_i, \lambda_i)$ are compact reduced labeled graphs representing the same GBS group $G$, for $i=1,2$. If $q(G) \cap \mathbb{Z}=1$, then $A_1$ and $A_2$ are related by slide moves.
\end{proposition}

\subsection{Index of a segment}Let $X$ be a locally finite $G$–tree. A $segment$ is an
edge path $\sigma = (e_1,\cdots,e_k)$ with no backtracking. 
Its initial and terminal vertices are $\partial_0(\sigma)=\partial_0(e_1)$ 
and $\partial_1(\sigma)=\partial_1(e_k)$,
respectively. The pointwise stabilizer of $\sigma$ is $G_\sigma =
G_{\partial_0(\sigma)} \cap G_{\partial_1(\sigma)}$. The index of $\sigma$ is the number $i(\sigma) = [G_{\partial_0 \sigma}  : G_\sigma]$. One can compute the index of any segment by applying the remark below iteratively, which can be found in \cite{forester2022incommensurable}.

\begin{remark}\label{index of segment}
     When $\sigma=(e_1, e_2)$ with $n_j=\lambda(e_j)$ and $m_j=\lambda(\overline{e}_j)$, for $j=1,2$; then $i(\sigma) = n_1 n_2/ \gcd(m_1,n_2)$. 
\end{remark}

\subsection{Subgroups of GBS groups}
If $G$ is a GBS group represented by a labeled graph $(A, \lambda)$, then there is a one-to-one
correspondence between conjugacy classes of GBS subgroups of $G$ (excluding hyperbolic cyclic subgroups) and admissible branched coverings $(B, \mu) \rightarrow (A, \lambda)$.
An $admissible$ $branched$ $covering$ from labeled graph $(B, \mu)$ to $(A, \lambda)$ consists of a surjective graph morphism 
$\pi: B \shortrightarrow A$ and a degree map 
$d: V(B) \sqcup E(B) \shortrightarrow \mathbb{N}$ satisfying $d(e)=d(\overline{e})$ for $e \in E(B)$, and if $e \in E(A)$ with $v=\partial_0(e)$, $u \in \pi^{-1}(v)$,  $k_{u,e}=\gcd(d(u), \lambda(e))$ then  
\begin{enumerate}
    \item $|\pi^{-1}(e) \cap E_0(u)|=k_{u,e}$
    \item If $e' \in \pi^{-1}(e) \cap E_0(u)$ then $\mu(e')=\lambda(e)/k_{u,e}$, and $d(e')=d(u)/k_{u,e}$.
\end{enumerate}
\begin{figure}[!ht]
\hspace*{3.5cm}
\begin{tikzpicture}
\small
  \draw[very thick] (1,4) -- (4,5);
  \draw[very thick] (1,4) -- (4,3);
  \draw[rotate around={-10:(1,4)},very thick] (1,4) -- (4,5);
  \filldraw[fill=white,thick] (1,4) circle (.7mm);
  \draw (0.8,4) node[anchor=east] {$u$};
  \draw[rotate around={-5:(1,4)}] (4.05,4.112) node[transform
  shape,anchor=center] {$\vdots$}; 
  \draw (4,4) node[anchor=west] {$\begin{rcases}  & \\ & \\ &
      \\ & \end{rcases} \ k_{u,e} \ = \
    \gcd(\textcolor{mygreen}{d(u)},\textcolor{violet}{\mu(e)})$}; 
  \draw[rotate around={18.435:(1,4)}] (1.15,4) node[transform
    shape,anchor=south west,violet] {$\mu(e)\big\slash k_{u,e}$};

  \draw[thick,->,color=mygreen] (1,3.3) -- (1,1.7);
  \draw[color=mygreen] (1,2.5) node[anchor=west] {$d(u)$}; 
  \draw[thick,->,color=mygreen] (2.6,2.83) -- (2.6,1.7);
  \draw[color=mygreen] (2.6,2.2) node[anchor=west] {$d(u)\big\slash
    k_{u,e}$};  
  
  \draw[very thick] (1,1) -- (4,1); 
  \filldraw[fill=white,thick] (1,1) circle (.7mm);
  \draw (0.8,1) node[anchor=east] {$v$};
  \draw (2.6,0.9) node[anchor=north] {$e$};
  \draw (1.15,1) node[anchor=south west,violet] {$\mu(e)$};
\end{tikzpicture}
\caption{The admissibility condition. Each edge of $p^{-1}(e) \cap
  E_0(u)$ has label $\mu(e)/k_{u,e}$ and degree $d(u)/k_{u,e}$. There
  are $k_{u,e}$ such edges.}\label{admissiblefig} 
\end{figure}

Let  $(A, \lambda)$ be a labeled graph, and $\pi: B \to A$ be a covering map in the topological sense. For any $e\in E(B)$, we can define the label $\mu(e)=\lambda(\pi(e))$. Then $(B, \mu)$ is labeled graph, and any constant degree map $c$, coprime to all edge labels of $A$, makes $\pi:(B, \mu) \to (A, \lambda)$ into an admissible cover.

\begin{remark}\label{ABC and top covering}
    An admissible branched cover from labeled graph $(B, \mu)$ to labeled graph $(A, \lambda)$ describes a topological covering of fibered $2$-complexes $Z_{(B, \mu) } \longrightarrow Z_{(A, \lambda)}$.
\end{remark}

Finite index subgroups of GBS groups without proper $p$--plateau have a nice description. The notion of a \textit{plateau} was introduced in \cite{levitt1}.
For a labeled graph $(A, \lambda)$ and a prime number $p$, a non-empty connected subgraph $P \subseteq A$ is a \emph{$p$--plateau} if for every edge $e \in E(A)$ with $v=\partial_0(e)$ belonging to the vertex set of $P$: $p \mid \lambda(e)$ if and only if $e \not \in E(P)$. The subgraph $P \subset A$ is called a \textit{plateau} if it is a $p$--plateau for some prime $p$. A plateau $P$ is considered \emph{proper} if $P\neq A$. 
 
\begin{proposition} \cite{levitt1}
\label{Levitt}
    Given a connected labeled graph $(A, \lambda)$, the following conditions are equivalent:
    \begin{itemize}
        \item every admissible covering $\pi: \overline{A} \to A$ is a topological covering;
        \item  $A$ contains no proper plateau.
    \end{itemize}
\end{proposition}

\subsection{Depth Profile}
The notion of depth profile is introduced in \cite{forester2022incommensurable} as a commensurability invariant of GBS groups. One defines an equivalence relation on the set of subsets of $\mathbb{N}$ by declaring that
$S \subset \mathbb{N}$ is equivalent to the set $\{n/ \gcd(r,n) : n \in S\}$  for each $r \in \mathbb{N}$ and taking the symmetric and
transitive closure. Let us denote the set $\{n/ \gcd(r,n) : n \in S\}$ by $S/r$.
\begin{proposition}\cite{forester2022incommensurable} \label{S_sim_S'}
    Two subsets $S$ and $S'$ are equivalent if and only if there exist $r, r' \in \mathbb{N}$ such that $S/r=S'/r'$.
\end{proposition}
Given $S \subset \mathbb{N}$ and $k\in \mathbb{N}$ such that $\gcd(x,k)=1$ for all $x \in S$, define 
\begin{equation}\label{S[k]}
S[k] = \{xk^i: x \in S, \text{ and } i \geq 0\}.
\end{equation}

\begin{lemma}\label{S(k)_equi_S'(k)}
Suppose $S, S' \subset \mathbb{N}$, $k \in \mathbb{N}$, and $\gcd(s,k)=\gcd(s',k)=1$ for all $s \in S$ and $s' \in S'$. Then $S$ and $S'$ are equivalent if and only if $S[k]$ and $S'[k]$ are equivalent.
\end{lemma}
\begin{proof}
Proposition \ref{S_sim_S'}  implies that if $S$ is equivalent 
to $S'$, then $S/r=S'/r'$ for some $r, r' \in \mathbb{N}$. We can assume that 
 $\gcd(r,k)=\gcd(r',k)=1$ since replacing $r$ by $r/\gcd(r,k)$ and $r'$ by $r'/\gcd(r',k)$ does not change the sets $S/r$ and $S'/r'$, as $k$ is coprime to all integers in $S$ and $S'$. This implies that $\gcd(r, sk^i)=\gcd(r,s)$ and $\gcd(r',s'k^i)=\gcd(r',s')$. Thus, the following sequence of equalities, together with Proposition \ref{S_sim_S'}, implies that $S[k]$ is equivalent to $S'[k]$. 

$$(S[k])/r=(S/r)[k]=(S'/r')[k]=(S'[k])/r'$$

For the converse, suppose $S[k]$ is equivalent to $S'[k]$. Then, by Proposition \ref{S_sim_S'}, we have $(S[k])/r=(S'[k])/r'$ for some $r, r' \in \mathbb{N}$.
We claim that $S/r=S'/r'$. Assuming the claim, $S$ is equivalent to $S'$ by Proposition \ref{S_sim_S'}. Since $S/r \subset (S[k])/r=(S'[k])/r'$, for every $s \in S$, we have
$$\frac{s}{\gcd(r,s)}=\frac{s'k^j}{\gcd(r',s'k^j)}=\frac{s'}{\gcd(r',s')} \frac{k^j}{\gcd(r',k^j)}$$
for some $j \geq 0$ and $s' \in S'$. This is possible only if $k^j/\gcd(r',k^j)=1$ since $\gcd(s,k)=1$. Therefore, $s/\gcd(s,r) \in S'/r'$. Thus, we have $S/r \subseteq S'/r'$, and similarly $S'/r' \subseteq S/r$. 
 \end{proof}

Let $G$ be a non-elementary GBS group, and $V < G$ is a non-trivial elliptic subgroup. For an element $g \in G$,  define its $V$--\textit{depth} as $D_g(V)=[V: V \cap V^g]$, where $V^g$ denotes the subgroup $\{gxg^{-1}: x \in V\}$. Define the \emph{depth profile} 
$$\mathscr{D}(G,V)=\{D_g(V): g \in G\ \text{is hyperbolic and } q(g)=\pm 1\} \subset \mathbb{N}.$$
The depth profile is commensurability invariant, i.e., if two non-elementary GBS groups  $G_1$ and $G_2$ are commensurable, then the sets $\mathscr{D}(G_1, V_1)$ and  $\mathscr{D}(G_2, V_2)$ are equivalent in the above sense (proved in \cite{forester2022incommensurable}).

\section{The $p$-modular homomorphism}
For a prime number $p$, the $p$--\textit{adic valuation} on the field of rational numbers is the map $\nu_p: \mathbb{Q}^* \to \mathbb{Z}$ defined by $\nu_p(\frac{a}{b})=\nu_p(a)-\nu_p(b)$, where $\nu_p(m)=\max\{e\in \mathbb{N}: p^e \mid m\}$.

\begin{definition} [\textbf{$\mathbf{p}$--modular homomorphism}]For a GBS group G and a prime number $p$, the $p$--\textit{modular homomorphism} is the map $q_p: G \to  \mathbb{Z}$ defined as $q_p=\nu_p \circ q$, where $q$ is the modular homomorphism on $G$, and $\nu_p$ is the restriction of the $p$-adic valuation on $Q^*$. The GBS group $G$ is called $p$--\textit{unimodular
} if $q_p(G)=0$.
\end{definition}

Let $(A, \lambda)$ be a labeled graph. For an edge path $(e_1, e_2, \cdots, e_k)$, define $q_p^A (e_1, e_2, \cdots, e_k) = \nu_p \circ q_A (e_1, e_2, \cdots, e_k)$, where we define $q_A$ by the right hand side in formula (\ref{eq0}).

\begin{lemma}\label{q_p^A for p-unimodular grp}
Let $G$ be a $p$-unimodular GBS group represented by labeled graph $(A, \lambda)$. For edge paths $(e_1, e_2, \cdots, e_k)$ and $(f_1, f_2, \cdots, f_l)$ with the same initial and terminal vertices, i.e., $\partial_0(e_1)=\partial_0(f_1)$ and $\partial_1(e_k)=\partial_1(f_l)$, we have
$$q^{A}_p (e_1, e_2, \cdots, e_k)= q^{A}_p (f_1, f_2, \cdots, f_k).$$
\end{lemma}

Therefore, for a fixed vertex $w_0 \in V(A)$, we get a well-defined function $h_p:V(A) \to \mathbb{Z}$ (with respect to $w_0$) defined by 

\begin{equation}\label{height function}
    h_p(v) = q_p^A (e_1, e_2, \cdots, e_k)
\end{equation} for any path $(e_1, e_2, \cdots, e_k)$ in $A$ from $w_0$ to $v$.

\begin{proof}
    Since $(e_1, e_2, \cdots, e_k)$ and $(f_1, f_2, \cdots, f_l)$ are two edge paths in $A$ from $w_0$ to $v$, the path $(e_1, \cdots, e_k, \overline{f}_l,\overline{f}_{l-1}, \cdots \overline{f}_1 )$ represents a $1-$cycle in $H_1(A)$. Then,
    \begin{align*}
       0 & = q_p(e_1, \cdots, e_k, \overline{f}_l,\overline{f}_{l-1}, \cdots \overline{f}_1 )\\
        &= q_p^A(e_1, \cdots, e_k)+q_p^A(\overline{f}_l,\overline{f}_{l-1}, \cdots \overline{f}_1 )\\
        &= q_p^A(e_1, \cdots, e_k)-q_p^A (f_1, f_2, \cdots, f_l)
     \end{align*}
    Thus, $q_p^A(e_1, \cdots, e_k)=q_p^A (f_1, f_2, \cdots, f_l)$.
\end{proof}

\begin{remark}
   We can define the map $q_n^A: G \to  \mathbb{Z}$ for any positive integer $n\geq 1$, and Lemma \ref{q_p^A for p-unimodular grp} also holds for the map $q_n^A$.
\end{remark}

\section{p-unimodularity and coverings}\label{construction of admissible cover for p-uni}

In this section, we present a general result regarding an admissible branched cover of a $p$-unimodular labeled graph. This result will be utilized in the subsequent section to demonstrate that all uniform lattices in Aut$(X_{dm,dn})$ are commensurable when $m$ and $n$ have no divisor less than or equal to $d$.

\begin{theorem}\label{coprime p cover}
Suppose $(A,\lambda)$ is a finite labeled graph which is $p$-unimodular. Then there exists a finite admissible cover $(A_p, \lambda_p)$ such that $p \nmid \lambda_p(e)$ for all $e \in E(A_p)$.
\end{theorem}

\begin{proof}
    First, perform expansion moves on the edges of $(A, \lambda)$ to obtain a labeled graph $(A_1, \lambda_1)$ such that if $p|\lambda_1(e)$ for some $e\in E(A_1)$, then $\lambda_1(e)=p$ and $\lambda_1(\overline{e})=1$. 
    
    If $e \in E(A)$ and $p \mid \lambda(e)$, then perform $\nu_p(\lambda(e))$    expansion moves as in Figure \ref{exp_(p,1)}. 
    \\

\begin{figure}[h]
    \centering
    \includegraphics[width = 0.75 \textwidth]{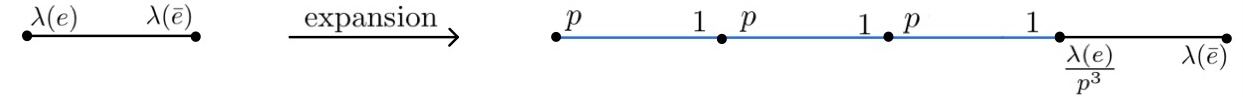}
    \caption{Example illustrating $\nu_p(\lambda(e))=3$ expansion moves for an edge $e \in E(A)$.}
    \label{exp_(p,1)}
\end{figure}

Consider the  subgraph $B_1$ of $A_1$ with the following vertex and edge set
\begin{itemize}
    \item $V(B_1)=V(A_1)$
    \item $E(B_1)=\{ e \in E(A_1) : p \nmid \lambda_1(e), \lambda_1(\overline{e})\}$
\end{itemize}
Note that if $e \in E(A_1)$ but $e \not \in E(B_1)$, then $\lambda_1(e)=p$ and  $\lambda_1(\overline{e})=1$ or $\lambda_1(e)=1$ and  $\lambda_1(\overline{e})=p$

Now we will define an admissible branched cover $(\Tilde{A}_1, \Tilde{\lambda}_1)$ of $(A_1, \lambda_1)$ with the property that $p \nmid \Tilde{\lambda}_1(e)$ for all $e\in E(\Tilde{A_1})$. Figure \ref{construction of covers} illustrates this construction with an example.  If $B_1=A_1$ (equivalently, $E(A_1)=E(B_1)$) then take $(\Tilde{A_1}, \Tilde{\lambda_1})=(A_1, \lambda_1)$.                                                              
\\

Claim: The function $h_p: V(A_1) \to \mathbb{Z}$ defined by equation \eqref{height function} with respect to a fixed vertex $w_0\in V(A_1)$ is constant on each component of $B_1$.

If $v_1, v_2 \in V(B_1)$ are in the same component of $B_1$ then there exists an edge path $(e_1, e_2, \cdots e_k)$ in $B_1$ with $\partial_0(e_1)=v_1$, $\partial_1(e_k)=v_2$. Since $e_i \in E(B_1)$,  $p \nmid \lambda(e) , \lambda(\overline{e})$. Therefore,
\begin{equation}
    q^{A}_p(e_1, e_1, \cdots e_k)=\nu_p \left( \frac{\lambda(e_1), \lambda(e_2), \cdots, \lambda(e_k)}{\lambda(\overline{e}_1), \lambda(\overline{e}_2), \cdots, \lambda(\overline{e}_k)}\right)=0
\end{equation}

Choose any edge path $(f_1, f_2, \cdots f_r)$ in $A_1$ from $w_0$ to $v_1$. Then the edge path given by $(f_1, f_2, \cdots f_r,e_1, e_1, \cdots e_k)$ has  initial vertex $w_0$  and terminal vertex $v_2$. Then, 
\begin{align*}
    h_p(v_2)
    &= q^{A}_p(f_1, f_2, \cdots f_r,e_1, e_1, \cdots e_k)\\
    &=q^{A}_p(f_1, f_2, \cdots f_r)+ q^{A}_p(e_1, e_1, \cdots e_k)\\
    &=q^{A}_p(f_1, f_2, \cdots f_r)\\
    &=h_p(v_1).
\end{align*}

This completes the proof of the claim. 
For $k \in \mathbb{Z}$, define $W_k$ to be the subgraph of $B_1$ spanned by vertices $v \in V(B_1)$ with $h_p(v)=k$. Since $B_1$ is a compact graph, and the $p$-modulus of each edge in $B_1$ is either $0$, $1$ or $-1$;  $W_k \neq \emptyset$  only for $k \in \{\alpha, \alpha+1 \cdots \alpha+\beta \}$ for some $\alpha, \beta \in \mathbb{Z}$. 

\begin{figure}[h]
    \centering
    \includegraphics[width = 0.75 \textwidth]{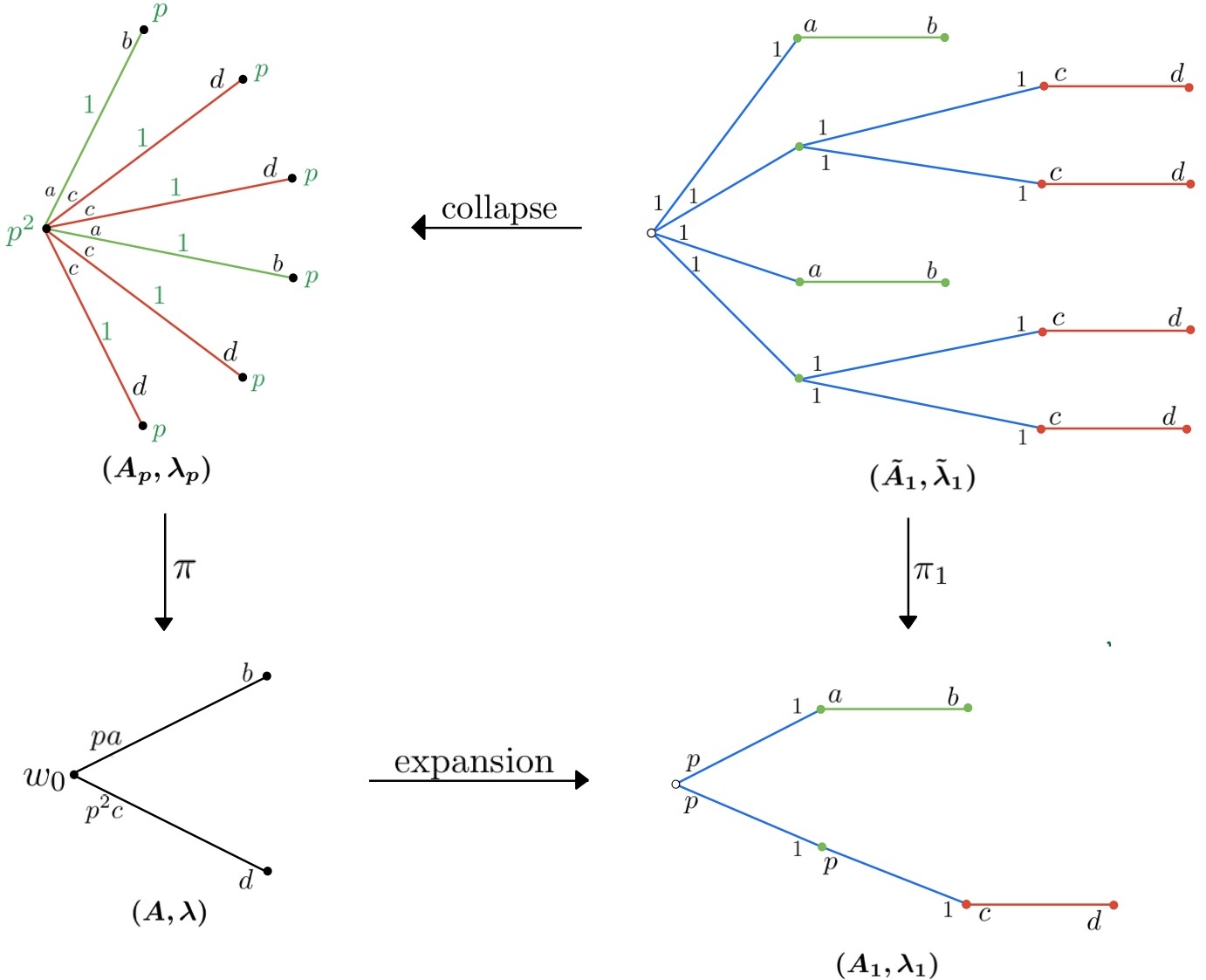}
    \caption{Construction of the admissible branched cover of the labeled graph $(A, \lambda)$ when $p=2$. Here $W_0$ is the white vertex, $W_1$ is the green subgraph and $W_2$ is the red subgraph. Hence $(\Tilde{A}_1, \Tilde{\lambda}_1)$ contains $2^0$ copies of $W_0$, $2^1$ copies of $W_1$, and $2^2$ copies of $W_2$.}
    \label{construction of covers}
\end{figure}

Construct an admissible branched covering $(\Tilde{A}_1,\Tilde{\lambda}_1 )$ of $(A_1, \lambda_1)$  by taking a disjoint union of $p^i$  copies of $W_{\alpha+i}$ for all $ 0\leq i \leq \beta$ together with some new edges. These copies are denoted as  $W_{\alpha+i}^1, W_{\alpha+i}^2, \cdots, W_{\alpha+i}^{p^i}$ with the same edge labels as $W_{\alpha+i}$. The surjective graph homomorphism  $\pi_1:\Tilde{A}_1 \to A_1$ for this admissible branched covering maps each copy $W_{\alpha+i}^j$ to $W_{\alpha+1}$ via the identity map. Define the degree  $d_1$ of each vertex and edge in $W_{\alpha+i}^j$ to be  $p^{\beta-i}$. Define
\\
$$V(\Tilde{A}_1)=\{v^j_{\alpha+i}: v\in V(W_{\alpha+i}), \text{ } 0\leq i \leq \beta,  \text{ and } 1\leq j \leq p^i\} .$$

Consider the surjective homomorphism 
$$\Phi: \frac{\mathbb{Z}}{p^{i+1}\mathbb{Z}} \to\frac{\mathbb{Z}}{p^i\mathbb{Z}}$$
defined by $\Phi([a])=[a]$ for all $[a]\in \frac{\mathbb{Z}}{p^{i+1}\mathbb{Z}}$.
The new edges in $\Tilde{A_1}$ are given  as follows: for each $e \in E(A_1)-E(B_1)$ with $\partial_0(e) \in V(W_{\alpha+i})$,  $\partial_1(e) \in V(W_{\alpha+i+1})$ for some $0 \leq i \leq \beta-1$, we have $\lambda(e)=p$,
and $\lambda(\overline{e})=1$. Then there are $p^{i+1}$ new edges in $\Tilde{A_1}$, $\{e^j: 1 \leq j \leq p^{i+1}\}$ with $\partial_1(e^j)=(\partial_1(e))^j_{\alpha+i+1}$ and $\partial_0(e^j)=(\partial_0(e))^{\phi(j)}_{\alpha+i}$. Define the label of these edges and their involution to be $1$. $\pi_1$ maps the edge  $e^j$ to $e$ with a degree $p^{\beta-i-1}$.

 Let $v \in V(W_{\alpha+i}) \subset V(A_1)$. Then,  for any $\Tilde{v} \in \pi_1^{-1}(v)$, we have $d_1(\Tilde{v})=p^{\beta-i}$. For an edge $e \in E(A_1)$ with $\partial_0(e)=v$ and $\partial_1(e)=w$, one of the following conditions holds: either $w \in  V(W_{\alpha+i})$, or $w \in  V(W_{\alpha+i+1})$, or $w \in  V(W_{\alpha+i-1})$. 
 
If $w \in  V(W_{\alpha+i})$, then $e \in E(W_{\alpha+i})$ and $ p\nmid \lambda_1(e)$, hence $\gcd(\lambda_1(e), d_1(v))=1$, and any $\Tilde{e} \in \pi_1^{-1}(e)$ is an edge in $W_{\alpha+i}^j$ for some $1 \leq j \leq p^i$. Therefore $|\pi_1^{-1}(e) \cap E_0(\Tilde{v})|=1$, $\Tilde{\lambda}_1(\Tilde{e})=\lambda_1(e)$, and $d_1(\Tilde{e})=p^{\beta-i}=d(\Tilde{v})$. 

If $w \in  V(W_{\alpha+i+1})$, then $\lambda_1(e)=p$, hence $\gcd(\lambda_1(e), d_1(v))=p$, and any $\Tilde{e} \in \pi_1^{-1}(e)$ is a new edge. Therefore $|\pi_1^{-1}(e) \cap E_0(\Tilde{v})|=p$, $\Tilde{\lambda}_1(\Tilde{e})=1$, and $d_1(\Tilde{e})=p^{\beta-i-1}$.

If $w \in  V(W_{\alpha+i-1})$, then $ \lambda_1(e)=1$, hence $\gcd(\lambda_1(e), d_1(v))=1$ and any $\Tilde{e} \in \pi_1^{-1}(e)$ is again a new edge. Therefore $|\pi_1^{-1}(e) \cap E_0(\Tilde{v})|=1$, $\Tilde{\lambda}_1(\Tilde{e})=1=\lambda_1(e)$, and $d_1(\Tilde{e})=p^{\beta-i}$. 
Hence $(\Tilde{A}_1, \Tilde{\lambda}_1)$ is an admissible branched cover of $(A_1, \lambda_1)$ for the surjective map $\pi_1:\Tilde{A}_1 \
\to A_1$, with the degree function $d_1$ defined above.

Since all new edges have label $1$, we can collapse all new edges to obtain a labeled graph $(A_p, \lambda_p)$. Also, $p \nmid \lambda_1(e)$ for all $e \in E(\Tilde{A_1})$, and $A_p$ is obtained by from $\Tilde{A}_1$ via collapse moves therefore $p \nmid \lambda_p(e)$ for all $e \in E(A_p)$.

Now, we will demonstrate that $(A_p, \lambda_p)$ defines an admissible branched cover of $(A, \lambda)$. We have the natural map $\pi:A_p \to A$ obtained from the following sequence of maps:
$$A_p \xrightarrow[]{expansion} \Tilde{A_1} \xrightarrow[]{\pi_1}  A_1 \xrightarrow[]{collapse} A$$

Define the degree of a vertex $u\in V(A_p)$ to be the maximum degree of all the vertices in $\Tilde{A}_1$ which maps to $u$ under the collapse map from $\Tilde{A}_1$ to $A_p$. Every edge $f \in E(A_p)$ corresponds to a unique edge $\Tilde{f} \in E(\Tilde{A}_1)$; define the degree of $f$ to be the degree of $\Tilde{f}$. Let $e \in A$, with $\partial_0(e)=v$, and $u \in \pi^{-1}(v)$. By the construction of $(\Tilde{A}_1, \Tilde{\lambda}_1)$, we have $\nu_p(d(u)) \geq \nu_p(\lambda(e))$, implying that $\gcd(\lambda(e), d(u))=p^{\nu_p(\lambda(e))}$. Furthermore, it is evident that $|\pi^{-1}(e) \cap E_0(u)|=p^{\nu_p(\lambda(e))}$, and for $e_p \in \pi^{-1}(e)\cap E_0(u)$,  we have $d(e_ p)=d(u)/p^{\nu_p(\lambda(e))}$, and $\lambda_p(e_p)=\lambda(e)/p^{\nu_p(\lambda(e))}$. Therefore the labeled graph $(A_p, \lambda_p)$ defines an admissible branched cover of $(A, \lambda)$ whose edge labels are coprime to $p$. This finishes the proof of the Theorem.

\end{proof}

\section{The Leighton property of $X_{m,n}$}\label{commensurability using leighton's property}
Assume $m$ and $n$ have no divisor less than or equal to $d$. We will show that for such a pair of numbers $(dm,dn)$, all torsion-free uniform lattices in Aut$(X_{dm,dn})$ are commensurable.

The main result of this section is Proposition \ref{mainprop_case3}, which implies that any torsion-free uniform lattice in Aut$(X_{dm,dn})$ for $(m,n)\neq (1,1)$ has a finite index subgroup represented by a directed labeled graph with edges labeled $m$ at the initial vertex and $n$ at the terminal vertex. 
It follows from Theorem \ref{prop4.6}(1) that the vertices of these graphs have $d$ incoming and $d$ outgoing edges incident to them. 
According to Leighton's theorem for graphs, any two labeled graphs with these properties share a common compact admissible cover.

Let's recall the statement of Leighton’s theorem here,
\begin{theorem}\cite[Leighton’s theorem]{leighton}
     Let $G_1$ and $G_2$ be finite connected graphs with a common cover. Then they have a common finite cover.
\end{theorem}

For the rest of the section, fix a general torsion-free uniform lattice $G$ in Aut$(X_{dm,dn})$, for $(m,n) \neq (1,1)$ unless otherwise stated. Let $(A, \lambda)$ be a compact GBS structure for $G$ given by Theorem \ref{prop4.6} with directed graph structure $E^+(A)$ and length function $l: E(A) \sqcup V(A) \longrightarrow \mathbb{N}$; then $A$ is strongly connected by Remark \ref{A is strongly connected}. 

\begin{lemma}\label{stru_of_A}
    If $e \in E^+(A)$, then $|\lambda(e)|=\alpha m$ and $|\lambda({\overline{e}})|=\beta n$ for some $1\leq \alpha, \beta \leq d$.
\end{lemma}
\begin{proof}
 Let $\partial_0(e)=v_1$ and $\partial_1(e)=v_2$. Since $A$ is strongly connected, there exists a directed cycle $(e_1, e_2, \cdots e_k)$ in $(A,\lambda)$, with $e_i \in E^+(A)$ for each $i$ and $e_1=e$. For $1 \leq i \leq k$, assume $\partial_0(e_i)=v_i$, and $\partial_1(e_i)=v_{i+1}$.
 
Claim(1): For any $i \in \{1, 2, \cdots k\}$, $l(v_i)
=\left(\frac{n}{m}\right)^{i-1} \left| \frac{\lambda(e_1) \lambda(e_2) \cdots\lambda(e_{i-1})}{\lambda(\overline{e}_1) \lambda(\overline{e}_2) \cdots\lambda(\overline{e}_{i-1})} \right| l(v_1)$.

We will prove the claim using induction on $i$. It is trivially true for $i=1$. Now,  assume that the claim holds for $i-1$. By Theorem \ref{prop4.6}(2),
$$l(\partial_0(e_{i-1})) |\lambda(e_{i-1})|= l(v_{i-1}) |\lambda(e_{i-1})|
=m l(e_{i-1}).$$
Therefore $l(e_{i-1})
=\frac{1}{m}l(v_{i-1})|\lambda(e_{i-1})|$. Furthermore, 
$$l(\partial_1(e_{i-1}) |\lambda(\overline{e}_{i-1}))|
=l(v_i) |\lambda(\overline{e}_{i-1})|
=n l(e_{i-1}).$$ 

Hence,

\begin{equation*}
\begin{split}
l(v_i)&=n \frac{l(e_{i-1})}{|\lambda(\overline{e}_{i-1})|}\\
&=\frac{n}{m}\frac{|\lambda(e_{i-1})|}{|\lambda(\overline{e}_{i-1})|}l(v_{i-1})\\
&=\frac{n}{m} \frac{|\lambda(e_{i-1})|}{|\lambda(\overline{e}_{i-1})|} \left(\frac{n}{m}
\right)^{i-2}\left|\frac{\lambda(e_1) \lambda(e_2) \cdots\lambda(e_{i-2})}{\lambda(\overline{e}_) \lambda(\overline{e}_2) \cdots\lambda(\overline{e}_{i-2}}) \right| l(v_1)\\
&=\left(\frac{n}{m}\right)^{i-1}\left|\frac{\lambda(e_1) \lambda(e_2) \cdots\lambda(e_{i-1})}{\lambda(\overline{e}_1) \lambda(\overline{e}_2) \cdots\lambda(\overline{e}_{i-1})} \right|l(v_1).
\end{split}
\end{equation*}
This proves the claim.
Since $(e_1, e_2, \cdots e_k)$ is a cycle, $v_1=v_{k+1}$. Hence,

\begin{align}
l(v_1)&=l(v_{k+1}) \nonumber\\
l(v_1)&= \left(\frac{n}{m}\right)^{k}\left|\frac{\lambda(e_1) \lambda(e_2) \cdots\lambda(e_{k})}{\lambda(\overline{e}_1) \lambda(\overline{e}_2) \cdots\lambda(\overline{e}_{k})} \right|l(v_1) \nonumber, 
\end{align}
and so
\begin{equation} \label{eq1}
     \left|\frac{\lambda(e_1) \lambda(e_2) \cdots\lambda(e_{k})}{\lambda(\overline{e}_1) \lambda(\overline{e}_2) \cdots\lambda(\overline{e}_{k})} \right| =\left( \frac{m}{n} \right)^{k}.
\end{equation}

Claim(2): $m \mid |\lambda(e_1)|$ and $n \mid |\lambda(\overline{e}_1)|$. 

Assuming claim(2) is true we get $|\lambda(e)|=\alpha m$, and  $|\lambda(\overline{e}_1)|=\beta n$, for some $\alpha, \beta \geq 1$. Also, by Proposition \ref{prop4.6}(1) $|\lambda(e)| \leq dm$, and $|\lambda(\overline{e}_1)| \leq dn$. Therefore,  $\alpha m=|\lambda(e)| \leq dm$, and $\beta n=|\lambda(\overline{e})| \leq dn$, implying  $\alpha, \beta\leq d$. This proves the lemma. 

To prove claim(2) we will only prove $m \mid |\lambda(e_1)|$, and $n \mid |\lambda(\overline{e}_1)|$ can be proved in similar manner. Let $m=p_1^{r_1}p_2^{r_2}\cdots p_k^{r_k}$ be the prime factorization of $m$. If $\nu_{p_i}(\lambda(e_1))\geq r_i$ for all $1 \leq i \leq k$ then $m \mid |\lambda(e)|$. Assume there is some $1 \leq i \leq k$ for which $\nu_{p_i}(\lambda(e_1)) < r_i$.

Without loss of generality we can assume $\nu_{p_1}(\lambda(e_1)) < r_1$. Then by equation \eqref{eq1}, $\nu_{p_1}(\lambda(e_i)) >r_1$ for some $2\leq i \leq k$. Again without loss of generality, we can assume that $\nu_{p_1}(\lambda(e_2)) >r_1$. Since $|\lambda(e_2)| \leq dm < p_1^{r_1+1}p_2^{r_2}\cdots p_k^{r_k}$, there exists $2 \leq j \leq l$ such that $\nu_{p_j}(\lambda(e_2)) < r_j$. Let $s_j=\nu_{p_j}(\lambda(e_2))$, then the set $J=\{j: s_j <r_j\}$
is a nonempty set. Proposition \ref{prop4.6}(2) provides the following sequence of implications;
$$l(v_2)|\lambda(e_2)|=ml(e_2)$$
$$\nu_{p_j}(l(v_2))+\nu_{p_j}(|\lambda(e_2)|)=\nu_{p_j}(m)+\nu_{p_j}(l(e_2))$$
$$\nu_{p_j}(l(v_2))=r_j-s_j+\nu_{p_j}(l(e_2))\geq r_j-s_j. $$ 
Therefore \begin{equation} \label{eq2}
k(v_2)=\gcd(l(v_2), m) \geq \prod_{j \in J}p_j^{r_j-s_j}.
\end{equation}

By Proposition \ref{prop4.6}(3), $E^+_0(v_2)=E^+_1 \sqcup \cdots \sqcup E^+_{k(v_2)}$ such that $\sum_{e \in E^+_i}|\lambda(e)|=\sum_{e \in E^+_j}|\lambda(e)|=C$ for all $1 \leq i,j \leq k_{v_2}$.
\begin{align*}
    dm&=\sum_{e \in E^+_0(v_2)}|\lambda(e)|\\
    &=\sum_{e \in E^+_1}|\lambda(e)| + \cdots + \sum_{e \in E^+_{k(v_2)}}|\lambda(e)|\\
    &=Ck(v_2)
\end{align*}
The fact that $e_2 \in E^+_0(v_2)$ and equation \ref{eq2} imply that,
\begin{equation}\label{eq3}
    |\lambda(e_2)| \leq C =\frac{dm}{k(v_2)} \leq \frac{dm}{\prod_{j\in J}p_j^{r_j-s_j}}=d\prod_{j\notin J}p_j^{r_j} \prod_{j\in J}p_j^{s_j}
\end{equation}

By our assumption, $|\lambda(e_2)|=cp_1^{r_1+\epsilon_1}p_2^{s_2}p_3^{s_3}\cdots p_k^{s_k}$ for some $c, \epsilon_1 \geq 1$. The fact that no prime less than or equal to $d$ divides $m$ gives the following inequality which contradicts inequality \ref{eq3}.
\begin{align*}
|\lambda(e_2)|&>dp_1^{r_1}p_2^{s_2}p_3^{s_3}\cdots p_k^{s_k}\\
&= dp_1^{r_1} \prod_{\substack{j \notin J\\j\neq1} }p_j^{s_j} \prod_{j\in J}p_j^{s_j}\\
&\geq d \prod_{\substack{j \notin J} }p_j^{r_j} \prod_{j\in J}p_j^{s_j}\\
\end{align*}

where the last inequality follows from the fact that $s_j \geq r_j$ for $j \notin J$.
\end{proof}

\begin{proposition}\label{prop_case3}
    The GBS group $G$ represented by $(A, \lambda)$ is $p$-unimodular for all primes $p \leq d$.
\end{proposition}
\begin{proof}
    Suppose $(e_1, e_2, \cdots e_k)$ is a cycle in $H_1(A)$ (it need not to be a directed cycle). Let
    $$\{e_1, e_2, \cdots e_k\}= 
  \{e_{i_1},e_{i_2}, \cdots, e_{i_r}\} \sqcup \{e_{j_1},e_{j_2}, \cdots, e_{j_s} \}$$
  where, $\{e_{i_1},e_{i_2}, \cdots, e_{i_r}, \overline{e}_{j_1},\overline{e}_{j_2}, \cdots, \overline{e}_{j_s} \} \subseteq E^+(A)$.  Let $\partial_0(e_i)=v_i$ and $\partial_1(e_i)=v_{i+1}$. Then Theorem \ref{prop4.6}(2), and the fact that $l(e)=l(\overline{e})$ imply that

\[
    l(v_{i+1})= 
\begin{cases}
    \frac{n}{m}
    \frac{\left|\lambda(e_i)\right|}{\left| \lambda(\overline{e}_i)\right|}l(v_i),& \text{if } e_i \in E^+(A)\\
    \frac{m}{n}  \frac{\left|\lambda(e_i)\right|}{\left| \lambda(\overline{e}_i)\right|}l(v_i),              & \text{if  } e_i \in E^-(A).
\end{cases}
\]
Since $(e_1, e_2, \cdots e_k)$ is a cycle, $v_1=v_{k+1}$ and the same proof as in claim(1) of Lemma \ref{stru_of_A} implies
$$l(v_1)=l(v_{k+1})=\left( \frac{n}{m} \right)^{r-s} \frac{\left|\lambda(e_{i_1})\lambda(e_{i_2}) \cdots \lambda(e_{i_r})\right|}{\left| \lambda(\overline{e}_{i_1})\lambda(\overline{e}_{i_2}) \cdots \lambda(\overline{e}_{i_r})\right|}   \frac{\left|\lambda(e_{j_1})\lambda(e_{j_2}) \cdots \lambda(e_{j_s})\right|}{\left| \lambda(\overline{e}_{j_1})\lambda(\overline{e}_{j_2}) \cdots \lambda(\overline{e}_{j_s})\right|} l(v_1). $$
Therefore, 
\begin{align*}
q(e_1, e_2, \cdots e_{k}) 
&=\frac{\left|\lambda(e_1)\lambda(e_2) \cdots \lambda(e_k)\right|}{\left| \lambda(\overline{e}_1)\lambda(\overline{e}_2) \cdots \lambda(\overline{e}_k)\right|} \\
&=\frac{\left|\lambda(e_{i_1})\lambda(e_{i_2}) \cdots \lambda(e_{i_r})\right|}{\left| \lambda(\overline{e}_{i_1})\lambda(\overline{e}_{i_2}) \cdots \lambda(\overline{e}_{i_r})\right|}   \frac{\left|\lambda(e_{j_1})\lambda(e_{j_2}) \cdots \lambda(e_{j_s})\right|}{\left| \lambda(\overline{e}_{j_1})\lambda(\overline{e}_{j_2}) \cdots \lambda(\overline{e}_{j_s})\right|}\\
&=\left( \frac{m}{n}\right)^{r-s}.
\end{align*}

Since $p \leq d$, and $\gcd(m,i)=\gcd(n,i)=1$ for all $1\leq i\leq d$, we have $\gcd(m,p)=\gcd(n,p)=1$. Therefore $\nu_p\left( \frac{m}{n}\right)=0$, and
\begin{align*}
    q_p(e_1, e_2, \cdots e_{k}) &= \nu_p\circ q (e_1, e_2, \cdots e_{k})\\
    &= \nu_p\left(\left( \frac{m}{n}\right)^{r-s}\right)\\
    &=(r-s)\left(\nu_p\left(\frac{m}{n}\right) \right)\\
    &=0.
\end{align*}
\end{proof}

\begin{proposition}\label{mainprop_case3}
    For a fixed prime number $p \leq d$, there exists a finite-index GBS group $H_p \leq  G$ with compact directed $GBS$ structure $(A_p, \lambda_p)$ and a directed graph structure $E(A_p)=E^+(A_p) \sqcup E^-(A_p)$ satisfying (1)-(3) of Theorem \ref{prop4.6}, such that for all $e\in E^+(A_p)$,
    \begin{align}\label{im, jn}
        |\lambda_p(e)|=\alpha_e m  \: \:\text{and} \:\: |\lambda_p(\overline{e})|=\beta_e n
    \end{align}
 for some $1\leq \alpha_e, \beta_e \leq d$, coprime to $p$.
\end{proposition}
\begin{proof}
     By Proposition \ref{prop_case3}, the graph $(A, \lambda)$ is $p$-unimodular. Therefore, Theorem \ref{coprime p cover} guarantees the existence of $H_p \leq G$. As $H_p$ itself is a uniform lattice in Aut$(X_{dm,dn})$, $(A_p, \lambda_p)$ admits a directed graph structure satisfying (1)-(3) of Theorem \ref{prop4.6}. Finally, Lemma \ref{stru_of_A} implies that the labels of edges in $A_p$ satisfy equation \ref{im, jn}.
\end{proof}

\begin{corollary}\label{coro1}
There exists a finite-index $GBS$ group $H \leq G$ with a compact  $GBS$ structure $(\Tilde{A}, \Tilde{\lambda})$, and a directed graph structure $E(\Tilde{A})=E^+(\Tilde{A}) \sqcup E^-(\Tilde{A})$ satisfying (1)-(3) of Theorem \ref{prop4.6}, such that for all $e\in E^+(\Tilde{A})$
    $$|\lambda(e)|=m,  \: \:\text{and} \:\: |\lambda(\overline{e})|=n.$$ 
\end{corollary}
\begin{proof}
    The result follows from applying Proposition \ref{prop_case3} iteratively for every prime $p\leq d$.
\end{proof}

\begin{corollary}\label{comm upto conjugacy}
    Any two torsion-free uniform lattices in Aut$(X_{dm,dn})$ are commensurable up to conjugacy.
\end{corollary}
\begin{proof}
    Let $G_1$ and $G_2$ be torsion-free uniform lattices in Aut$(X_{dm,dn})$. By Corollary \ref{coro1}, there exist finite index GBS groups $H_i \leq G_i$ with compact GBS structures $(\Tilde{A_i}, \Tilde{\lambda_i})$ for $i\in \{1,2\}$ satisfying
\begin{equation}\label{eq5}
    |\lambda(e)|=m \: \:\text{and} \:\: |\lambda(\overline{e})|=n
\end{equation} 
for all $e \in E^+(\Tilde{A}_i)$. By Theorem \ref{prop4.6}$(2)$ we also have
\begin{equation}\label{eq6}
    \sum_{e\in E^+_0(v)} |\lambda(e)|=dm, \:\:\text{and} \sum_{e\in E^-_0(v)} |\lambda(\overline{e})|=dn.
\end{equation}

These equations imply that $\Tilde{A_1}$ and $\Tilde{A_2}$ are $2d$ regular graphs. Consequently, by Leighton's graph covering theorem, they also share a common finite-sheeted topological cover. Equation (\ref{eq6}) guarantees that this common cover can be labeled to create a common admissible cover of $\Tilde{A_1} $ and $\Tilde{A_2}$. Thus, $H_1$ and $H_2$ 
 are commensurable up to conjugacy in Aut$(X_{dm,dn})$, implying the same for $G_1$ and $G_2$. 

 \end{proof}

\subsection{ Torsion-free uniform lattices in Aut$(X_{d,d})$} This subsection addresses the last component of Theorem \ref{main theorem}, namely that any two torsion-free uniform lattices in 
Aut$(X_{d,d})$ are commensurable.

If Aut$(X_{m,n})$ is a discrete group, it cannot contain incommensurable lattices \cite[1.7]{basslubotzky}. Also, Aut$(X_{m,n})$ is discrete if and only if $\gcd(m,n)=1$ \cite[Theorem 4.8]{forester2022incommensurable}.  Therefore,  Aut$(X_{1,1})$ is discrete, and   any two lattices in Aut$(X_{1,1})$ are commensurable.

Let $\Gamma$ be a torsion-free uniform lattice in Aut$(X_{d,d})$ for $d \geq 2$. Then, by Proposition \ref{freely and compactly}, $\Gamma$ acts on $X_{d,d}$ freely and cocompactly, providing a covering space action.  Each branching line in $X_{d,d}$ covers a circle, and each strip covers either an annulus or a M\"{o}bius band, therefore the quotient space may not be orientable. However, by [\cite{bass}, Proposition 6.3] we can find an index $2$ subgroup $G$ of $\Gamma$ that acts on $X_{d,d}$ without changing the sides of any strip.  
Therefore, the quotient obtained from the action of $G$ on $X_{d,d}$ is a fibered $2$-complex $Z_{(A.\lambda)}$ for some labeled graph $(A, \lambda)$. For $e \in E(A)$, if $l(e)$ denotes the number of $2$ cells tiling the annulus corresponding to $e$ and $l(v)$ denotes the combinatorial length of the circle corresponding to $v$, then we have the following result which is derived from the proof of Proposition 4.6 of \cite{forester2022incommensurable}.

\begin{proposition}\label{X_d,d -G is uni}
  Suppose $G$ is a torsion-free uniform lattice in Aut$(X_{d,d})$ for $d \geq 2$ such that the action of $G$ on $X_{d,d}$ does not change the sides of any strip. Let $(A, \lambda)$ be the labeled graph structure such that $X_{d,d}/G$ is homeomorphic to $Z_{(A, \lambda)}$ with the associated length function $l: V(A) \sqcup E(A) \mapsto \mathbb{N}$. Then, for any edge $e \in E(A)$ we have, 
  \begin{enumerate}
      \item $l(\partial_0(e)) |\lambda(e)|= l(\partial_1(e)) |\lambda(\overline{e})|$.
      \item $G$ is unimodular.
  \end{enumerate}
\end{proposition}

\begin{proof}
Condition (1) follows from the cell structure of the annulus corresponding to the edge $e$. This annulus is tiled by $(1,1)$ cells whose boundary curves have length $l(e)$ which wrap $|\lambda(e)|$ and $|\lambda(\overline{e})|$ times onto the circles corresponding to $\partial_0(e)$ and $\partial_1(e)$, respectively.
\\
   To prove (2),  let $(e_1. \cdots, e_k)$ be a cycle  in $A$. Then, by (1)
    \begin{align*}
        q(e_1. \cdots, e_k)&= \frac{|\lambda(e_1) \cdots \lambda(e_k)|}{|\lambda(\overline{e}_1) \cdots \lambda(\overline{e}_k)|}\\
        &=\frac{l(\partial_1
        (e_1)) \cdots l(\partial_1
        (e_1))}{l(\partial_0
        (\overline{e}_1)) \cdots l(\partial_0
        (\overline{e}_1))}\\
        &=1
    \end{align*}
    where the last equality fallows form the fact that $\partial_0(e_{i+1})=\partial_1(e_i)$ for $1 \leq i \leq k-1$ and $\partial_1(e_k)=\partial_0(e_1)$.

\end{proof}

\begin{proposition}\label{levitt automorphism}\cite{levitt-auto}
    If $G$ is a non-elementary GBS group then $G$ is unimodular if and only if it has a finite index subgroup isomorphic to $F
_n \times \mathbb{Z}$ for some $n>1$.

\begin{proposition}
    Any two torsion-free uniform lattices in Aut$(X_{d,d})$, for $d \geq 2$ are commensurable. 
\end{proposition}
\begin{proof}
    Let $\Gamma_1$ and $\Gamma_2$ are two torsion-free uniform lattices in Aut$(X_{d,d})$. We can assume that $\Gamma_i$ act on $T_{d,d}$ without inversion (possibly after passing to an index $2$ subgroup of $\Gamma_i$) for $i=1,2$.
    By Proposition \ref{X_d,d -G is uni} $\Gamma_i$
are unimodular, and hence contain a finite index subgroup isomorphic to $F_{n_i} \times \mathbb{Z}$ for some $n_i >1$ by Proposition \ref{levitt automorphism}. Finally, since $F_{n_1} \times \mathbb{Z}$  and $F_{n_2} \times \mathbb{Z}$ are commensurable, it follows that $\Gamma_1$ and $\Gamma_2$ are commensurable.

\end{proof}
\end{proposition}

\subsection{Leighton's Property} \begin{definition}
    We say that a cell complex $X$ has the $Leighton$ $property$ if every pair of compact cell complexes, both having $X$ as their common universal cover, admits a common finite-sheeted covering.
    
\end{definition}

By Leighton's theorem, trees have the Leighton property.

 \begin{theorem}\label{one direction of 1.4}
     When $m$ and $n$ have no prime divisor less than or equal to $d$, the cell complex $X_{dm,dn}$ satisfies the Leighton property.
 \end{theorem}
 \begin{proof}
     
  For cell complexes $X_1$ and $X_2$ with common universal cover $X_{dm,dn}$, the fundamental groups $\pi_1(X)$ and $\pi_2(X)$ are torsion-free since $X_1$ and $X_2$ are finite-dimensional aspherical cell complexes, hence defining torsion-free uniform lattices in Aut$(X_{dm,dn})$. For $i=1,2$, we can choose labeled graphs $(A_i, \lambda_i)$ with associated fibered 2-complexes $X_i$.

 For $(m,n) \neq (1,1)$, the results of Proposition \ref{mainprop_case3} and its corollaries \ref{coro1}, \ref{comm upto conjugacy} imply  the labeled graphs $(A_1, \lambda_1)$ and $(A_2, \lambda_2)$ admit a common finite sheeted admissible branched covering. The fibered $2$-complex associated with this common finite sheeted admissible branched covering provides a common finite sheeted cover of $X_1$ and $X_2$ by Remark \ref{ABC and top covering}.

For $(m,n) = (1,1)$ and $d=1$, all lattices in  Aut$(X_{1,1})$ are commensurable as it is a dicrete group. Therefore $X_{1,1}$ satisfies Leighton's property.

  For $(m,n) = (1,1)$ and $d \geq 2$, we can assume that $X_i$ are orientable (possible after passing to degree $2$ covering of $X_i$) for $i=1,2$. ${X}_i$ is a fibered $2$-complex associated to some labeled graph $(A_i, \lambda_i)$ which is unimodular by Proposion \ref{X_d,d -G is uni}, hence $p$-unimodular for all prime $p$. By Theorem \ref{coprime p cover},  we can find a $2d$ regular admissible branched cover $(\Tilde{A}_i,\Tilde{\lambda}_i )$ of $(A_i, \lambda_i)$ with edge labels $1$. By Leighton's graph covering theorem, the graphs $\Tilde{A}_1$ and $\Tilde{A}_2$ admit a common finite topological cover, denoted as $\Tilde{A}$. Assigning labels $1$ to each edge of $\Tilde{A}$ yields a labeled graph $(\Tilde{A}, \Tilde{\lambda})$ which defines an admissible branched cover of $(\Tilde{A}_i,\Tilde{\lambda}_i )$. Since the composition of admissible branched covers is again an admissible branched cover, $(\Tilde{A}, \Tilde{\lambda})$ defines a common admissible branched cover for $(A_1, \lambda_1)$ and $(A_2, \lambda_2)$. Thus, the fibered $2$-complex associated with $(\Tilde{A}, \Tilde{\lambda})$ provides a common finite sheeted topological cover of $X_i$.

\end{proof}

\section{Example of incommensurable lattices using the depth profile}\label{eg using DP}
We provide examples of incommensurable lattices in both Case (I) and Case (II), utilizing the commensurability invariant known as the depth profile.

\subsection{Incommensurable lattices in Case (I)}
   In this subsection, we will provide examples of incommensurable lattices in Aut$(X_{dm,dn})$, when there is a prime number $p \leq d$ such that $p \mid d$, and $p\mid m$ or $p\mid n$. Without loss of generality, let $p$ be a prime number which divides both $m$ and $d$. Note that $p \nmid n$ since $\gcd(m,n)=1$. 

\begin{figure}[h]
    \centering
    \includegraphics[width = 0.37 \textwidth]{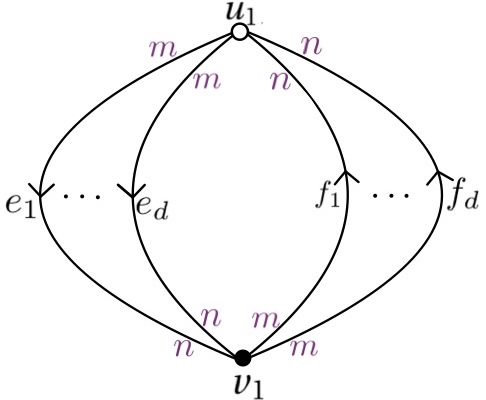}
    \caption{$B_1$}
    \label{fig:B_1}
\end{figure}

\textbf{The lattice {$\mathbf{\Gamma_1}$}.} 
Consider the lattice $\Gamma_1$ defined by the directed labeled graph  $(B_1, \mu_1)$ in Figure (\ref{fig:B_1}). It is a bipartite graph with two vertices $u_1$ (white) and $v_1$ (black), and $2d$ directed edges. The edges $e_1, e_2, \cdots e_d$ are directed from $u_1$ to $v_1$ and the edges $f_1, f_2, \cdots, f_d$ are directed from $v_1$ to $u_1$. We have $\mu_1(e_i)=\mu_1(f_i)=m$ and $\mu_1(\overline{e_i})=\mu_1(\overline{f_i})=n$.

\begin{figure}[h]
    \centering
    \includegraphics[width=0.4\linewidth]{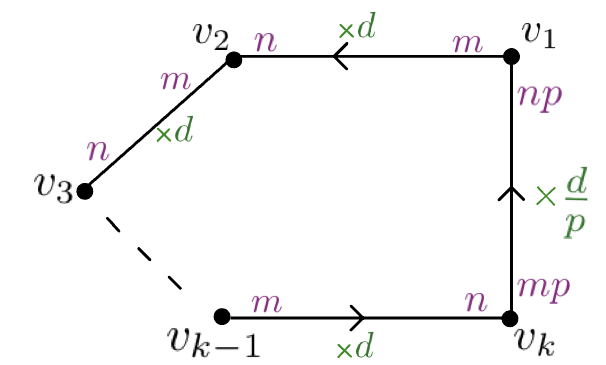}
    \caption{ Graph $B_l$ for $l\geq 2$, defining lattices in Aut$(X_{dm,dn})$ for a prime number $p \mid d,m$. Here, $\times \alpha$ denotes the number of edges between vertices.}
    \label{fig:B_2}
\end{figure}

\textbf{The lattice {$\mathbf{\Gamma_k}$} for {$\mathbf{k\geq 2}$}.} 
This group is defined by the directed labeled graph $(B_k, \mu_k)$ in Figure (\ref{fig:B_2}). It has $k$ vertices $v_1, \cdots , v_k$ , and $d(k-1)+ \frac{d}{p}$ directed edge. There are $d$ directed edges from $v_{i}$ to $v_{i+1}$ for $1 \leq i \leq k-1$ with initial label $m$ and terminal label $n$ and there are $\frac{d}{p}$ directed edges from $v_k $ to $v_1$ with initial label $pm$ and terminal label $pn$. \\

The groups $\Gamma_k$ for $k \geq 1 $ are all latties in Aut($X_{dm,dn}$) by theorem \ref{labled graphs indeed define lattices}. We will compute the depth profiles of $\Gamma_k$, for $k \geq 1$, and will show that these depth profiles are not equivalent using Lemma \ref{S(k)_equi_S'(k)}. This will prove that these groups are pairwise incommensurable.
\begin{definition}
    A segment $\sigma$ in a $G-$tree is called $unimodular$ if $i(\sigma)=i(\overline{\sigma})$.
\end{definition}

\begin{proposition} \cite{forester2022incommensurable}\label{index and depth profile}
    Let $G$ be a GBS group 
 and $X$ be its GBS tree. Suppose $V$ is the stabilizer of a vertex $x \in V(X)$. Define the set
 $$\mathscr{I}(x)=\{i(\sigma): \sigma \text{ is a non-trivial unimodular segment with endpoints in } Gx\}.$$
 
 Then 
 $$\mathscr{D}(G,V) \subseteq \mathscr{I}(x) \subseteq \mathscr{D}(G,V) \cup \{1\}.$$
\end{proposition}

\begin{proposition}

Let $X_1$ be the Bass-Serre trees for the given GBS structure on $\Gamma_1$. Suppose $V_1$ is the stabilizer of a vertex $x_1 \in V(X_1)$ which maps to the black vertex in the graph of groups $B_1$. Then $\mathscr{I} (x_1)=S_1[n]-\{1\}$ where 
    $$S_1=\big{\{} m^i: i \in \mathbb{N} \cup \{0\} \big{\}}.$$
   
\end{proposition}
\begin{proof}

Note that the unimodular segments in $X_1$ with both endpoints in $\Gamma_1x_1$ have even lengths, and the vertices along these segments alternate between black and white. Additionally, for this entire proof, we only consider the unimodular segments from a black to another black vertex in $X_1$.

\begin{figure}[h]
    \centering
    \includegraphics[width=0.4\linewidth]{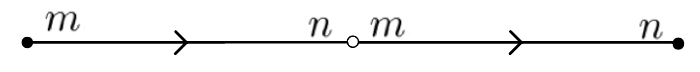}
   \caption{A segment $\tau$ in $X_1$.}
    \label{fig:tau_j}
\end{figure}

Let's denote any length $2$ segment in $X_1$ with the labels given in Figure \ref{fig:tau_j} by $\tau$. For such segments $\tau$, the initial index $i(\tau)$ is $m^2$ and the terminal index $i(\overline{\tau})$ is $n^2$. Now every unimodular segment $\sigma$ in $X$ of length $>2$ has one of the following forms:

    \begin{enumerate}[label={(\arabic*)}]
  \item $ \sigma= \sigma_1 \sigma_2$ with $\sigma_1$, $\sigma_2$ unimodular segments in $X$
  \item $\sigma= \tau \sigma_1 \overline{\tau}$ with  $\sigma_1$ a unimodular segment in $X$
  \item $\sigma= \overline{\tau} \sigma_1 {\tau}$ with  $\sigma_1$ a unimodular segment in $X$.
\end{enumerate}

Let $D_1$ denote the set $S_1[n]-\{1\}$ (see equation \eqref{S[k]} for definition of $S[k]$). It is easy to verify that $D_1$ is closed under taking $\lcm$.
We will show that every unimodular segment $\sigma$ has index in $D_1$ by induction on its length. 
Let's denote any edge in $X_1$ with an initial label $m$ and a terminal label $n$ by $e$. Then,
the length $2$ unimodular segments in $X_1$  whose end vertices are black are $e\overline{e}$ and $\overline{e}e$. By Remark \ref{index of segment} we can see that $i(e\overline{e})=m \in D_1$, and $i(\overline{e}e)=n \in D_1$.

 If $\sigma$ is of type(1), then by Remark \ref{index of segment}, we have $i(\sigma)=\lcm(i(\sigma_1), i(\sigma_2)) \in D_1$.
 
 By Remark \ref{index of segment},  and using the fact that $\gcd(\lcm(a, b), b)=b$, if $\sigma$ is of type(2), then $i(\sigma)=n^{-2}m^2 \lcm(i(\sigma_1),n^2) \in D_1$, and If $\sigma$ is of type(3), then $i(\sigma)=m^{-2}n^2\lcm(i(\sigma), m^2) \in D_1$. This shows that $\mathscr{I}(x_1) \subseteq D_1$.
 
Finally, consider the following unimodular segments in $X_1$:
 \begin{enumerate}[label={(\alph*)}]
 \item $e^i \overline{e}^i$
 \item $\overline{e}^j e^j $

  \end{enumerate}
 
These segments have indices $m^{i}$, and $n^{j}$, respectively. Therefore the index of a concatenation of segments of type(a) and type(b) is $\lcm (m^i, n^j)=m^i n^j$. Hence we also have $D_1 \subseteq \mathscr{I}(x_1)$. 
\end{proof}

\begin{proposition}
    Let $X_k$ be the Bass-Serre tree for the given GBS structure for $\Gamma_k$, for $k \geq 2$. Suppose $V_k$ is the stabilizer of a vertex $x_k$ which maps $v_1$. Then $\mathscr{I}(x_k)=S_k[n]-\{p\}$ where  $$S_k=\{pm^{ik}, m^{ik+1}, \cdots  m^{ik+(k-1)} : i \geq 0\}$$
\end{proposition}
\begin{proof}
    Note that, the unimodular segments in $X_k$ with both endpoints in $\Gamma_k x_k$ have even lengths as they contain equal numbers of forward- and backward-oriented edges.

Let $D$ denote the set $S_k[n]-\{p\}$. It is easy to check that $D$ is closed under taking lcm. We will show that every unimodular segment in $X_k$  with both endpoints in $\Gamma_k x_k$ has index in $D$ by induction on its length. Let's denote edges in $X_k$ with initial label $m$ and terminal label $n$ by $e$, and edges with initial label $pm$ and terminal label $pn$ by $f$.

For the base case, we will show that the index of unimodular segments of length $2l$,  for $1 \leq l \leq k$, is contained in $D$. Observe that $i(\tau_1e\overline{e}e\overline{e}\tau_2)=i(\tau_1 e\overline{e}\tau_2)$ for any segments $\tau_1$ and $\tau_2$. Therefore it is sufficient to compute the index of unimodular segments that are either $e^l\overline{e}^l$, $\overline{f}\overline{e}^{l-1}e^{l-1}f$ or concatenations of smaller unimodular segments. Since $i(e^l\overline{e}^l)=m^l$ and $i(\overline{f}\overline{e}^{l-1}e^{l-1}f)=pn^l$, and since $D$ is closed under taking $lcm$, it follows from Remark \ref{index of segment} that the set of indices of unimodular segments of length $2l$ is contained in $D$.

\begin{figure}[h]
    \centering
    \includegraphics[width=0.6\linewidth]{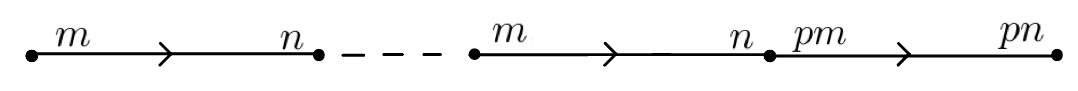}
    \caption{A segment $\tau$ in $X_k$}
    \label{segment tau}
\end{figure}

 Consider the segment $\tau=e^{k-1}f$ as in Figure \ref{segment tau}. Note that $i(\tau)=pm^k$ and  $i(\overline{\tau})=pn^k$.  The index of every unimodular segment in $X_k$ of length $>2k$ is contained in the index of one of the following unimodular segments;

    \begin{enumerate}[label={(\arabic*)}]
  \item $\sigma \sigma'$ with $\sigma$, $\sigma'$ unimodular segments in $X_k$
  \item $\tau \sigma \overline{\tau}$ with  $\sigma$ a unimodular segment in $X_k$
  \item $\overline{\tau} \sigma {\tau}$ with  $\sigma$ a unimodular segment in $X_k$.
\end{enumerate}
The index of an unimodular segment of type(1) is  $\lcm(i(\sigma), i(\sigma'))$ which is contained in  $D$ as it is closed under taking lcm. Also, by Remark \ref{index of segment} and the fact that $gcd(lcm(a, b), b) = b$, we get $i(\tau \sigma \overline{\tau})=(pn)^{-k}pm^k \lcm({pn^k, i(\sigma)}) \in D$, and $i(\overline{\tau} \sigma \tau)=(pm)^{-k}pn^k \lcm({pm^k, i(\sigma)}) \in D$. This shows that $\mathscr{I}(x_k) \subseteq D$.\\
Finally, consider the following unimodular segments in $X_k$ for $i \geq 0$, and $0 \leq j \leq k-1
$;
 \begin{enumerate}[label={(\alph*)}]
\item $\tau^i \overline{\tau}^i$
 \item $\tau^i e^j \overline{e}^j\overline{\tau}^i$
 \item $\overline{\tau}^i  e^j \overline{e}^j \tau^i$
 \item $\overline{\tau}^i \overline{f} \overline{e}^{j-1} e^{j-1}  f \tau^i$
  \end{enumerate}
  These segments have indices $pm^{ik}$, $m^{ik+j}$, $pn^{ik}$, and $pn^{ik+j}$, respectively. Therefore,
\begin{itemize}
    \item Concatenation of segments of type(a) and type(c) has index $pm^{i_1k}n^{i_2 k}$, 
    \item Concatenation of segments of type(a) and type(d) has index $pm^{i_1k}n^{i_2 k+j_2}$,
    \item Concatenation of segments of type(b) and type(c) has index $m^{i_1k+j_1}n^{i_2k}$,
    \item Concatenation of segments of type(b) and type(d) has index $m^{i_1k+j_1}n^{i_2k+j_2}$.
\end{itemize}
Hence we also have $D \subseteq \mathscr{I}(x_k)$. 
\end{proof}

\begin{corollary}\label{depth profile of Gamma1}
    $\mathscr{D}(\Gamma_1, V_1)= S_1[n]-\{1\}$ and 
    $\mathscr{D}(\Gamma_k, V_k)=S_k[n]-\{p\}$ for $k \geq 2$.
\end{corollary}

\begin{proof}
    Since the set $\mathscr{I}(x_k)$ does not contain $1$ for $k \geq 1$, by Proposition \ref{index and depth profile} we have $\mathscr{D}(\Gamma_k, V_k)=\mathscr{I}(x_k)$.
\end{proof}

\begin{corollary}\label{ratio not agreeing}
If the prime number $p$ divides both $m$ and $d$, then the lattices $\Gamma_k \in \text{Aut}(X_{dm,dn})$  are pairwise abstractly incommensurable for $k\geq 1$.
\end{corollary}

\begin{proof}
   Enumerate the elements of $S_i$ in order for $i \geq 1$ and notice that each element divides the next one. Taking the ratio of successive elements we obtain the sequences $(m,m,m, \cdots)$ for $S_1$ and  
   $$\Bigl( \frac{m}{p},
  \underbrace{
    m, m , \cdots, m,
    }_{\text{$k-2$~elements}}
     pm, \frac{m}{p},
     \underbrace{
    m, m, \cdots, m,
    }_{\text{$k-2$~elements}}
    pm, \frac{m}{p},
     \underbrace{
    m, m , \cdots, m,
    }_{\text{$k-2$~elements}}
    \cdots
  \Bigr)$$
for $S_k$, $k \geq 2$. The tails of these ratio sequences are unchanged when passing from $S_i$ to $S_i/r$ for any $r\in \mathbb{N}$ because the values of $\gcd(r, m^j)$, and $\gcd(r, pm^{kj})$ stabilize as 
   $j \shortrightarrow \infty$, all to the same number. The tail for $S_1$ will never agree with the tail for $S_k$ as $p \neq 1$, so we get $S_1$ not equivalent to $S_k$ for $k \geq 2$.  Also, for $k, l \geq 2$, the tails for $S_l$ and $S_k$ will agree if and only if  $l=k$. Using Lemma \ref{S_sim_S'}, we conclude that $S_l[n]$ is not equivalent to $S_k[n]$ for $k\neq l$.
  Furthermore, it's evident that $S_1[n]$ is equal to the set
   $(S_1[n]-\{1\})/n$, which is equivalent to $\mathscr{D} (\Gamma_1, V_1)=S_1[n]-\{1\}$. Similarly, $S_k[n]$ is equal to the set $(S_k[n]-\{p\})/n$, which is equivalent to $\mathscr{D} (\Gamma_k, V_k)=S_k[n]-\{p\}$. Therefore by Lemma \ref{S(k)_equi_S'(k)} the depth profiles of $\Gamma_l$ and $\Gamma_k$ are not equivalent for $k \neq l$ and hence these groups are not abstractly commensurable.
   
\end{proof}

\subsection{Incommensurable lattices in Case (II)}\label{new section}
In this section, we provide examples of lattices in Aut$(X_{dm,dn})$ that are abstractly incommensurable when $m$ or $n$ is $1$, and there exists  $p <d$ (not necessarily a prime), $m,n\neq p$ such that either $p \mid m$ or $p \mid n$. Without loss of generality, we can assume $ m = 1$, $n \neq p$, and $p \mid n$. \cite{forester2022incommensurable} provides an example of incommensurable lattices in Aut$(X_{d,dn})$ when $n \neq  p$, using depth profile as a commensurability invariant (see Theorem \ref{Forester's thm}(2)). Building on this, we will construct infinitely many such examples using the depth profile.

\begin{figure}[h]
    \centering
    \includegraphics[width=0.45\linewidth]{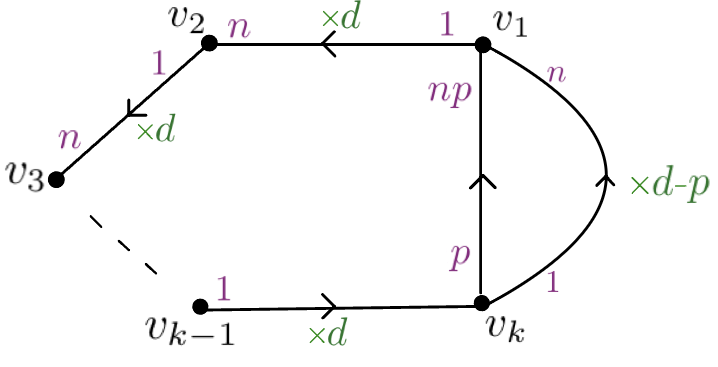}
    \caption{The graph $D_k$ for $k\geq 2$, defining a lattice in Aut$(X_{d,dn})$ for $p<d$.}
    \label{Delta_k}
\end{figure}

\textbf{The lattice {$\mathbf{\Delta_k}$} for {$\mathbf{k \geq 2}$} .}
Consider the group $\Delta_k$ defined by the directed labeled graph   $(D_k, \delta_k)$ shown in Figure(\ref{Delta_k}). The graph $D_k$ consists of $k$ vertices $v_1, v_2, \cdots v_k$ and 
$dk-p+1$ directed edges. There are $d$ directed edges from $v_i$ to $v_{i+1}$ for $1 \leq i \leq k-1$ and $d-p$ directed edges from $v_k$ to $v_1$, each with initial label $1$ and terminal label $n$. Additionally, there is a directed edge from $v_k$ to $v_1$ with initial label $p$ and terminal label $np$. 

Performing a sequence of collapse and slide moves on $D_k$ for $k \geq 2$ (see figure \ref{forth case pic} for an example with $k=5$), we obtain the following:

\begin{figure}[h]
    \centering
    \includegraphics[width=0.9\linewidth]{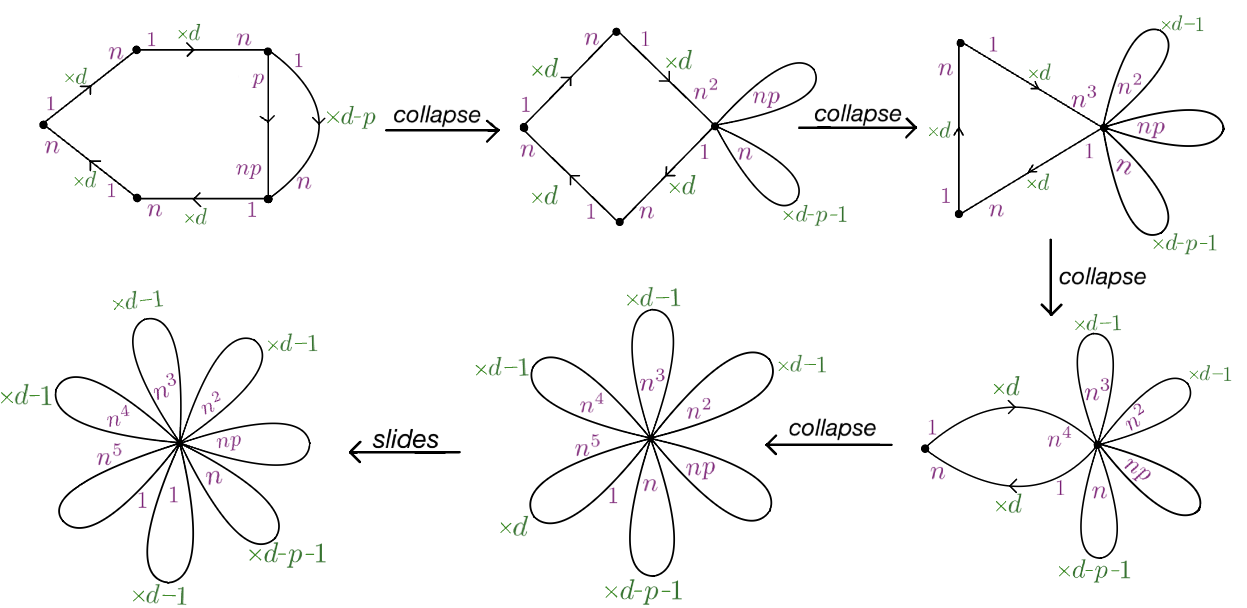}
    \caption{A sequence of collapse and slide moves on $(\Delta_5, \delta_5)$ that gives a bouquet of circles with the fundamental group $\Delta_5$. A number $z$ inside the
petal indicates that both ends of the petal are labeled $z$, while $\times y$ above a petal denotes
the multiplicity of that petal.}
    \label{forth case pic}
\end{figure}

\begin{align}
   \begin{split}
    \Delta_k &\cong   
   BS(pn,pn) \vee \bigvee_{d} BS(1, n^k) \vee \bigvee_{d-p-1 } BS(n,n) \vee \bigvee_{d-1} BS(n^2,n^2) \vee \bigvee_{d-1} BS(n^{3},n^{3}) \vee \cdots \\
     &\qquad  \bigvee_{d-1} BS(n^{k-1},n^{k-1}) 
   \end{split} \nonumber
   \\[1ex]
    \begin{split}
    &\cong BS(pn,pn) \vee BS(1, n^k) \vee \bigvee_{d-1 }BS(1,1) \vee   \bigvee_{d-p-1 } BS(n,n) \vee \bigvee_{d-1} BS(n^2,n^2) \vee \bigvee_{d-1} BS(n^{3},n^{3})   \\
   &\qquad \cdots \bigvee_{d-1} BS(n^{k-1},n^{k-1}). \label{bouquet instruction of delta_k}
    \end{split}
\end{align}

The following result from \cite{forester2022incommensurable} will be used to compute the depth profile of $\Delta_k$:

\begin{proposition} \cite{forester2022incommensurable}\label{depth profile of bouquets of circles}
    Let $G=BS(1,N) \vee \bigvee_{i=1}^r BS(n_i, n_i)$ for some $r \geq 1$, $N>1$, and $n_i$ dividing $N$. Suppose the set $\{n_1, n_2, \cdots n_r\}$ is closed under taking \emph{lcm} and contains $1$. Then, for the vertex group $V$,
    $$\mathscr{D}(G, V)= \{n_iN^j : j \geq 0 \textit{ and } 1 \leq i \leq r\}.$$
\end{proposition}

\begin{proposition}
    The set of groups $\{ \Gamma_1, \Delta_k: k \geq 2\}$ defines pairwise incommensurable lattices in \emph{Aut}$(X_{d,dn})$ when $n \geq d$, $gcd(n,d)=1$, and $n$ has a divisor $p<  d$.
 
\end{proposition}

\begin{proof}
   For $k\geq 2$, the group $\Delta_k$ defines a lattice in Aut$(X_{d,dn})$ by Theorem \ref{labled graphs indeed define lattices}.
   Let $V_k$ be the vertex group of the GBS structure \ref{bouquet instruction of delta_k} for $\Delta_k$. By Proposition \ref{depth profile of bouquets of circles}, the depth profile of $\Delta_k$ is:
   \begin{equation}
\mathscr{D}(\Delta_k, V_k)=
    \begin{cases}
\{pn^{ki+1}, n^{ki}, n^{ki+1}, n^{ki+2}, \cdots n^{ki+k-1} : i \geq 0\} & \textit{if } d > p+1\\
    
        \{pn^{ki+1}, n^{ki}, n^{ki+2}, \cdots n^{ki+k-1} : i \geq 0\} & \textit{if } d = p+1
    \end{cases}
   \end{equation}
   In both cases, enumerating the elements of $\mathscr{D}(\Delta_k, V_k)$ in ascending order and computing the successive ratios we get the periodic sequence:
  \begin{center}

   $\begin{cases}
       n, p, \frac{n}{p}, \underbrace{n , \cdots, n}_{k-1 \textit{ times}}, p, \frac{n}{p}, \underbrace{n , \cdots, n}_{k-1 \textit{ times}}\cdots & \textit{ if } d  >p+1\\

       pn, \frac{n}{p}, \underbrace{n , \cdots, n}_{k-2 \textit{ times}}, pn, \frac{n}{p}, \underbrace{n , \cdots, n}_{k-2 \textit{ times}}\cdots & \textit{ if } d =p+1.\\
   \end{cases}$
\end{center}
By Corollary \ref{depth profile of Gamma1} $$\mathscr{D}(\Gamma_1, V_1)=\{n^i: i\geq 0\},$$
and the ratio of successive terms in $\mathscr{D}(\Gamma_1, V_1)$ is $(n, n, n, \cdots )$. As the tail of these ratio sequences never agrees which is a commensurability invariant as mentioned in Corollary \ref{ratio not agreeing}, we conclude that any two groups in the set $\{\Gamma_1, \Delta_k: k \geq 2\}$ are incommensurable. 
\end{proof}

\section{Revisiting the CRKZ invariant}\label{eg using crkz}
In this section, we define a class of GBS groups and provide a necessary condition for two
groups in this class to be abstractly commensurable. This condition will enable us to prove that the set of groups $\{\Gamma_1, \Delta_k: k\geq 2\}$ defines pairwise incommensurable lattices in Aut$(X_{d,dn})$ when $n<d$ and $n \nmid d$. Furthermore, we also demonstrate that this condition is sufficient for a special subset of this class.

In \cite{MR4359918}, the authors introduced a class of GBS groups and 
constructed an isomorphism invariant for groups in this class (see subsection \ref{CRKZ_comm} for details). We will refer to this invariant as the CRKZ invariant.
Here we will give a new description of the CRKZ invariant for a larger class of GBS groups and prove the scaling property of the CRKZ invariant for finite index subgroups arising from topological covers (Theorem \ref{scaling property}). It follows that the CRKZ invariant is also a complete commensurability invariant.

Fix an integer $l \geq 1$. Suppose $G$ is a non-elementary GBS group whose image under the modular homomorphism $q: G \to \mathbb{Q}^*$ is generated by $1/n^{lL}$ for some $n \in \mathbb{N}$ and  $L \in \mathbb{Z}$. Suppose $G$ is represented by a labeled graph $(A, \lambda)$, and for all edges $e \in E(A)$, $\lambda(e)=n^{i_e}$ for $i_e \geq 0$.  To each GBS group in this form, we will associate a vector $\Vec{X}^l(G) \in (\mathbb{N} \cup \{0\})^l$ well defined up to cyclic permutation. We call this vector the \textit{length $l$ CRKZ invariant} of $G$.

\begin{definition}\label{level of vertex and edge}
Suppose $v_0 \in V(A)$ is a fixed vertex.
\begin{enumerate}
    \item 
    A vertex $v \in V(A)$ has $level$ $i$ with respect to the base vertex $v_0$ if for any path $(e_1,e_2, \cdots e_r) $ from $v_0$ to $v$, $q^A_n(e_1,e_2, \cdots e_r) \equiv i \pmod{l}$.
    \item  An edge $e \in E(A)$  has $level$ $i$ with respect to base vertex $v_0$ if, for any path $(f_1, f_2, \cdots, f_s) $ from $v_0$ to $\partial_0(e)$,  $\nu_{n}(\lambda(e))+q^A_n(f_1, f_2, \cdots, f_s) \equiv i \pmod{l}$.
    \end{enumerate}
\end{definition}

Let $V_{v_0}^i(A)$ (and $E_{v_0}^i(A)$) denote the set of vertices (and edges) that have level $i$ with respect to the base vertex $v_0$. In particular, 

$$V_{v_0}^i(A) =\{ v \in V(A) :  q^A_n(e_1,e_2, \cdots e_r) \equiv i\pmod{l}\}$$
$$E_{v_0}^i(A) = \{ e \in E(A) : \nu_n(\lambda(e))+q^A_n(f_1, f_2, \cdots, f_s) \equiv i \pmod{l}\}$$
where  $(e_1,e_2, \cdots ,e_r)$ is an edge path in $A$ from $v_0$ to $v$, and $(f_1, f_2, \cdots, f_s)$ is an edge path in $A$ from $v_0$ to $\partial_0(e)$. Note that
 $V_{v_0}^{i+kl}(A)=V_{v_0}^i(A)$ and $E_{v_0}^{i+kl}(A)=E_{v_0}^i(A)$ for all $1\leq i\leq l$ and $k \in \mathbb{N}$.\\

 Definitions \ref{level of vertex and edge} is independent of the choice of path from $v_0$ to $v$ since if $(e_1,e_2, \cdots ,e_r)$  and $(f_1, f_2, \cdots, f_s)$ are two paths in $A$ from $v_0$ to $v$, then $(e_1,e_2, \cdots ,e_r, \bar{f}_s, \bar{f}_{s-1}, \cdots, \bar{f}_1)$ is a $1$-cycle in $H_1(A)$.
 Therefore 
 \begin{align*}  
     q^{A}_n(e_1,e_2, \cdots ,e_r)- q^{A}_n(f_1, f_2, \cdots, f_s) &=
     q^{A}_n(e_1,e_2, \cdots ,e_r, \bar{f}_s, \bar{f}_{s-1} \cdots \bar{f}_1)\\
     &\equiv 0 \pmod{l}
 \end{align*}
implies $ q^{A}_n(e_1,e_2, \cdots ,e_r) \equiv q^{A}_n(f_1, f_2, \cdots, f_s) \pmod{l}$. \\
\\

\begin{remark}
    Edges $e$ and $\overline{e}$ always have the same level, so $e \in E^i_{v_0}(A)$ if and only if  $\overline{e} \in E^i_{v_0}(A)$. In particular, $|E^i_{v_0}(A)|$ is even.
\end{remark}

\begin{definition}

Let $(A, \lambda)$ be a labeled graph. For a fixed integer $l\geq 1$ and $v_0 \in V(A)$, define a vector $\Vec{X}^l_{v_0}(A)$ in $(\mathbb{N} \cap \{0\})^l$ as

$$\Vec{X}^l_{v_0}(A) =\left(|E^{0}_{v_0}(A)|-2|V^{0}_{v_0}(A)|,|E^{1}_{v_0}(A)|-2|V^{1}_{v_0}(A)|, \cdots, |E^{l-1}_{v_0}(A)|-2|V^{l-1}_{v_0}(A)|\right). $$
The next lemma shows that these vectors are independent of the choice base point, up to cyclic permutation.

\end{definition}

\begin{lemma}\label{X(A) is indepent of choice of base point}
Let $v_0$ and $ w_0$ be two vertices in $A$. Then the vector $\Vec{X}^l_{w_0}(A)$ is a cyclic permutation of $\Vec{X}^l_{v_0}(A)$. In particular, $\sigma^{i_0}\Vec{X}^l_{v_0}(A)=\Vec{X}^l_{w_0}(A)$, where $\sigma=(l, l-1, \cdots, 1)$ is the cyclic permutation of $\{1, 2, \cdots, l\}$ and $i_0=q^A_n(e_1,e_2, \cdots ,e_r)$ for an edge path $(e_1,e_2, \cdots ,e_r)$ from $v_0$ to $w_0$.

\end{lemma}

\begin{proof}
 For an edge path $(f_1, f_2, \cdots, f_s)$ from $w_0$ to $v$, $(e_1,e_2, \cdots ,e_r, f_1, f_2, \cdots, f_s)$  is a path from $v_0$ to $v$ with
\begin{align*}
    q^A_n(e_1,e_2, \cdots ,e_r, f_1, f_2, \cdots, f_s)&=q^A_n(e_1,e_2, \cdots ,e_r)+q^A_n(f_1, f_2, \cdots, f_s) \\
    &\equiv i_0 + q^A_n(f_1, f_2, \cdots, f_s) \pmod{l}
\end{align*}
Therefore, we have  $V_{w_0}^i(A)=V_{v_0}^{i+i_0}(A)$ and $E_{w_0}^i(A)=E_{v_0}^{i+i_0}(A)$. Together with the facts $|V_{v_0}^{i+kl}(A)|=|V_{v_0}^i(A)|$ and $|E_{v_0}^{i+kl}(A)|=|E_{v_0}^i(A)|$ for all $1\leq i\leq l$ and $k \in \mathbb{Z}$ we get,

\begin{align*}
\sigma^{i_0} \left(\Vec{X}^l_{v_0}(A) \right)
&=\sigma^{i_0} \left(|E^{0}_{v_0}(A)|-2|V^{0}_{v_0}(A)|,|E^{1}_{v_0}(A)|-2|V^{1}_{v_0}(A)|, \cdots, |E^{l-1}_{v_0}(A)|-2|V^{l-1}_{v_0}(A)|\right)\\
&=\Big( |E^{i_0}_{v_0}(A)-2|V^{i_0}_{v_0}(A)|,|E^{i_0+1}_{v_0}(A)|-2|V^{i_0+1}_{v_0}(A)|, \cdots, \\  & \;\;\;\;\;\;\; |E^{i_0+l-1}_{v_0}(A)|-2|V^{i_0+l-1}_{v_0}(A)|  \Big)
\\
&=\left(|E^0_{w_0}(A)|-2|V^{0}_{w_0}(A)|,|E^1_{w_0}(A)|-2|V^{1}_{w_0}(A)|, \cdots, |E^{l-1}_{w_0}(A)|-2|V^{l-1}_{w_0}(A)| \right)\\
&=\Vec{X}^l_{w_0}(A).
\end{align*}
\end{proof}

 We will denote any element in the subset $\{\Vec{X}^l_{v_0}(A): v_0 \in V(A)\} \subset {(\mathbb{N}\cup \{0\})}^l$ as $\Vec{X}^l(A)$. In the view of Lemma \ref{X(A) is indepent of choice of base point},  $\Vec{X}^l(A)$ is a well-defined vector in ${(\mathbb{N}\cup \{0\})}^l$ up to cyclic permutation.

 Any two splittings of a non-elementary GBS group are in the same deformation space. This means that these graphs of groups are related via a sequence of expansion and collapse moves. Thus, by the next lemma, for a non-elementary GBS group  $G$,  we can associate a vector $\Vec{X}^l(G)$ well defined up to cyclic permutation. Sometime we will refer $\Vec{X}^l(G)$ as $\Vec{X}^l(A)$.

\begin{lemma}
    If two labeled graphs $A$ and $A'$ are in the same deformation space then $\Vec{X}^l(A)$ is a cyclic permutation of $\Vec{X}^l(A')$.
\end{lemma}
\begin{proof}
    It suffices to show that if $A'$ is obtained from $A$ via a collapse move, then $\Vec{X}^l(A')$ is a cyclic permutation of $\Vec{X}^l(A)$. 
\begin{figure}[h]
    \centering
    \includegraphics[width = 0.6 \textwidth]{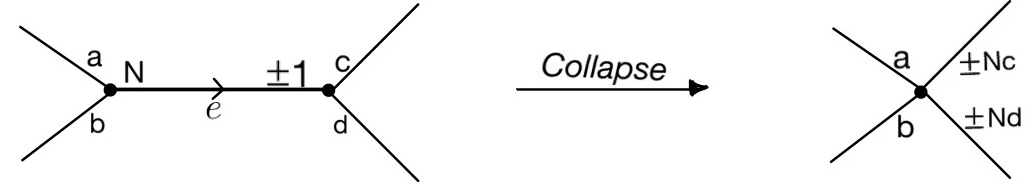}
    
\end{figure}

Suppose $A'$ is obtained from $A$ by collapsing an edge $e$ with $\lambda(e)=N$, and $\lambda(\bar{e})=1$. Let us denote the image of a vertex $v \in A$ via collapse move by $v' \in A'$. Choose $\partial_0(e)$ and $(\partial_0(e))'$
  to be the base vertex in $A$ and $A'$, respectively.   
  If $\nu_n(N) \equiv i_0 \pmod{l}$, then the vertex $\partial_1(e)$ has level $i_0$, and edges $e, \bar{e}$ also have level $i_0$. Therefore 

\[
     \left|V^i(A')\right|= 
\begin{cases}
    \left|V^i(A)\right|-1,              & \text{if } i \equiv i_0 \pmod{l}\\
    
    \left|V^i(A)\right|
   ,& \text{if } i \not \equiv i_0 \pmod{l},
\end{cases}
\]
and
\[
     \left|E^i(A')\right|= 
\begin{cases}
    \left|E^i(A)\right|-2,              & \text{if } i \equiv i_0 \pmod{l}\\
    
    \left|E^i(A)\right|
   ,& \text{if } i \not \equiv i_0 \pmod{l}.
\end{cases}
\]

Hence $|E^{i_0}(A')|-2|V^{i_0}(A')|=|E^{i_0}(A)|-2-2(|V^{i_0}(A)|-1)=|E^{i_0}(A)|-2|V^{i_0}(A)|$, and $\Vec{X}^l(A')=\Vec{X}^l(A)$ follows from the definition.

   \end{proof}

\begin{theorem}
    \label{scaling property}
    Let $G$ be a GBS group defined by labeled graph $(A, \lambda)$ without a proper plateau such that for all $e \in E(A)$, $\lambda(e)=n^{i_e}$ for some $i_e \geq 0$ and $q(G)= \langle (1/n^{lL})\rangle_{\mathbb{Q}^*}$. Then for every index $d$ subgroup $H \leq G$, $\Vec{X}^l(H)$ is a cyclic permutation of $d\Vec{X}^l(G)$.
\end{theorem}

\begin{proof}
    Since $A$ does not contain a proper $p$--plateau, the GBS group $H$ is represented by a labeled graph $\Tilde{A}$ for some $d$-sheeted topological covering $\pi: \Tilde{A} \to A$, by Proposition \ref{Levitt}. Fix a base vertex $v_0 \in A$ and $\Tilde{v_0}\in \pi^{-1}(v_0) \subseteq V(\Tilde{A})$.
    
   For any $\Tilde{v} \in V(\Tilde{A})$  and a directed edge path $(\Tilde{e}_1,\Tilde{e}_2, \cdots, \Tilde{e}_r)$ from $\Tilde{v_0}$ to $\Tilde{v}$, $(\pi(\Tilde{e}_1),\pi(\Tilde{e}_2), \cdots, \pi(\Tilde{e}_r))$
    is directed edge path in $A$ from $v_0$ to $\pi(\Tilde{v})$ with $q_n^{\Tilde{A}}(\Tilde{e}_1,\Tilde{e}_2, \cdots, \Tilde{e}_r)=q_n^A(\pi(\Tilde{e}_1),\pi(\Tilde{e}_2), \cdots, \pi(\Tilde{e}_r))$.
   Therefore, for $\Tilde{v} \in \pi^{-1}(v)$ and $v \in V^i(A)$, we have $\Tilde{v} \in V^i(\Tilde{A})$.
   Similarly, for $\Tilde{e} \in \pi^{-1}(e)$ and $e \in E^i(A)$, we have $\Tilde{e} \in E^i(\Tilde{A})$.
   Thus, $|E^i(\Tilde{A})|=d|E^i(A)|$ and $|V^i(\Tilde{A})|=d|V^i(A)|$. Now $\Vec{X}^l_{\Tilde{v}_0}(\Tilde{A})=d\Vec{X}^l_{v_0}(A)$ follows from the definition of $\Vec{X}^l$. 
\end{proof}

\begin{corollary}\label{dV(A)}
    Suppose $G_i$ are the GBS group represented by the directed labeled graphs $(A_i, \lambda_i)$ without proper $p$--plateau, for $i=1,2$. If $G_1$ is commensurable to $G_2$, then $c_1\Vec{X}^l(G_1)$ is a cyclic permutation of $c_2\Vec{X}^l(G_2)$ for some $c_1, c_2 \in \mathbb{N}$.
\end{corollary}
\begin{proof}
    Suppose $G_1$ is commensurable to $G_2$, then for some finite index subgroup $H_i \leq G_i$, $H_1$ is isomorphic to $H_2$. Let $H_i$ be represented by the labeled graph $B_i$ with the covering map $\pi_i :B_i \to A_i$. Then $B_1$ and $B_2$ are related via a sequence of expansion and collapse moves. Therefore, $\Vec{X}^l(B_1)$ is a cyclic permutation of $\Vec{X}^l(B_2)$.
   Also, $[G_i:H_i] \Vec{X}^l(G_i)=\Vec{X}^l(B_i)=\Vec{X}^l(H_i)$ by  Proposition \ref{dV(A)}. Thus, if $G_2$ is commensurable to $G_2$, then $c_1\Vec{X}^l(G_1)$ is a cyclic permutation of $c_2\Vec{X}^l(G_2)$ for $c_i=[G_i:H_i]$.
\end{proof}

\subsection{Incommensurable lattices in Case (III)} 
Next, we prove that any two lattices in $\{\Gamma_1, \Delta_k: k\geq 2\}$ are incommensurable when $p=n$.

Recall the lattice $\Gamma_1$ defined by Figure \ref{fig:B_1} for $m=1$. By
collapsing the edge $e_d$, and performing $d-1$ slide moves, we find that 
\begin{align*}
    \Gamma_1 &\cong \bigvee_{d} BS(1, n^2) \vee \bigvee_{d-1} BS(n,n)\\
    &\cong BS(1, n^2) \vee \bigvee_{d-1 }BS(1,1) \vee \bigvee_{d-1} BS(n,n) 
\end{align*}

Also, the group $\Delta_k$  for $n=p$ from \eqref{bouquet instruction of delta_k} is

\begin{align*}
    \begin{split}
    \Delta_k &\cong BS(n^2,n^2) \vee BS(1, n^k) \vee \bigvee_{d-1 }BS(1,1) \vee   \bigvee_{d-n-1 } BS(n,n) \vee \bigvee_{d-1} BS(n^2,n^2) \vee \bigvee_{d-1} BS(n^{3},n^{3})   \\
   &\qquad \cdots \bigvee_{d-1} BS(n^{k-1},n^{k-1}). \label{bouquet instruction of delta_k}
    \end{split}
\end{align*}

\begin{figure}[h]
\centering
  \subcaptionbox{$(B_1',\mu_1')$\label{fig3:a}}{\includegraphics[width=1.7in]{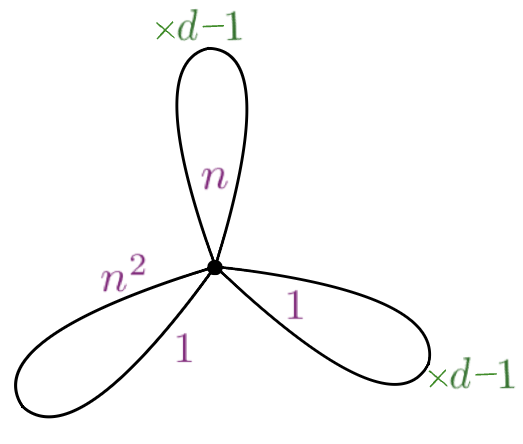}}\hspace{1em}%
  \subcaptionbox{$(D_2',\delta_2')$\label{fig3:b}}{\includegraphics[width=1.7in]{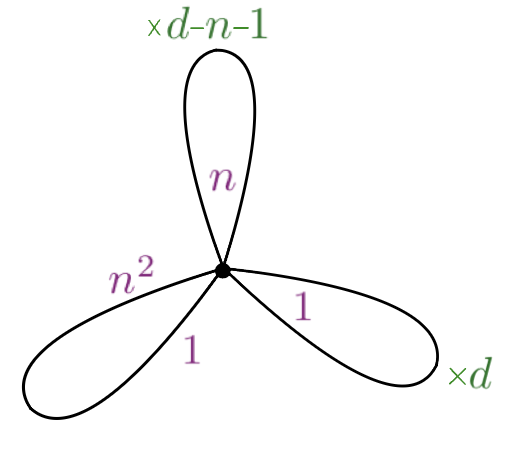}}\hspace{.5em}%
  \subcaptionbox{$(D_k',\delta_k')$, for $k\geq 3$\label{fig3:c}}{\includegraphics[width=2.2in]{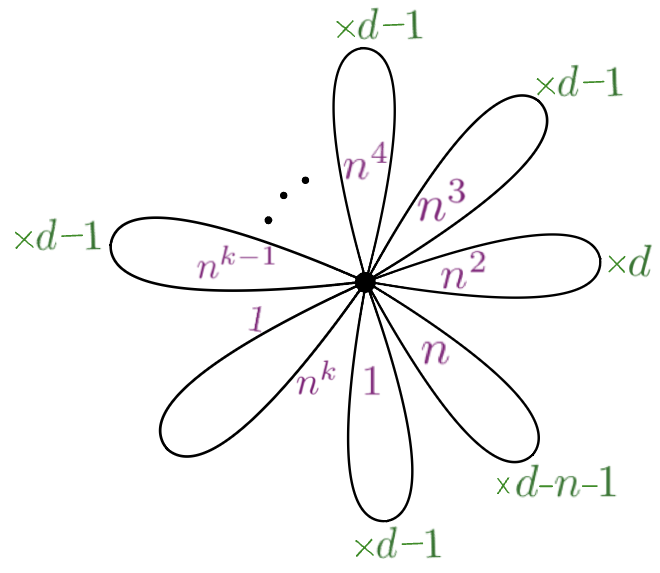}}
   \caption{Labeled graphs representing the groups $\Gamma_1$ and $\Delta_k$ for $k \geq 2$.}%
  \label{bouquet for Gammas}%
\end{figure}

Therefore, $\Gamma_1$ and $\Delta_k$ are the fundamental groups of the labeled graphs which are bouquets of circles. These bouquets consist of one loop labeled $1$ and $n^k$, along with additional loops labeled $n^i$ for $0 \leq i \leq k-1$, as shown in Figure \ref{bouquet for Gammas}.

\begin{proposition}
    The set of groups $\{ \Gamma_1, \Delta_k : k \geq 2 \}$ define pairwise incommensurable torsion-free uniform lattices in \emph{Aut}$(X_{d,dn})$, when $n < d$ and $\gcd(n,d)=1$.
\end{proposition}
\begin{proof}
   To show that any two groups in $\{ \Gamma_1, \Delta_k: k \geq 2 \}$ are incommensurable, we will demonstrate that they contain incommensurable finite index subgroups.

To obtain an index $l$ subgroup $\Delta_{k,l}$ in $\Delta_k$, unwind the loop in the graph $(D_k', \delta_k')$ labeled $1$ and $n^k$ into a circle of  length $l$ (see Figure \ref{index l subgroup of Gamma} for an example). Next, collapse all but one of the edges labeled $1$ and $n^k$ to obtain a bouquet of circles with edges labeled $n^i$ for $0 \leq i \leq kl$. The fundamental group of this labeled graph is $\Delta_{k,l}$. Similarly, unwinding the loop labeled $1$ and $n^2$ in $(B_1', \mu_1')$ into the circle of length $l$  gives index $l$ subgroup in $\Gamma_{1,l}$ in $\Gamma_1$. It is straightforward to see that the modular homomorphism of $\Delta_{k,l}$ is generated by $\frac{1}{n^{kl}}$, and the vector $\Vec{X}^{kl}(\Delta_{k,l})$ is as follows:
\begin{enumerate}

 \item   $\Vec{X}^{2l}(\Delta_{2,l})=2( 
 d, d-n-1, d, d-n-1, \cdots, d, d-n-1 )\in \mathbb{N}^{2l}$ 
  
  \item    $\Vec{X}^{kl}(\Delta_{k,l})=2( d-1, d-n-1, d, \underbrace{ d-1, \cdots, d-1}_{k-3 \text{ elements}}, \cdots
 d-1, d-n-1, d, \underbrace{ d-1, \cdots, d-1}_{k-3 \text{ elements}})$ in $\mathbb{N}^{kl}$ for $k \geq 3$.
\end{enumerate}

Meanwhile, the modular homomorphism of $\Gamma_{1,l}$ is generated by $\frac{1}{n^{2l}}$, and 
$$\Vec{X}^{2l}(\Delta_{1,l})=2( 
 d-1, d-1, \cdots, d-1, d-1) \in \mathbb{N}^{2l}.$$

\begin{figure}[h]
    \centering
    \includegraphics[width=.9\linewidth]{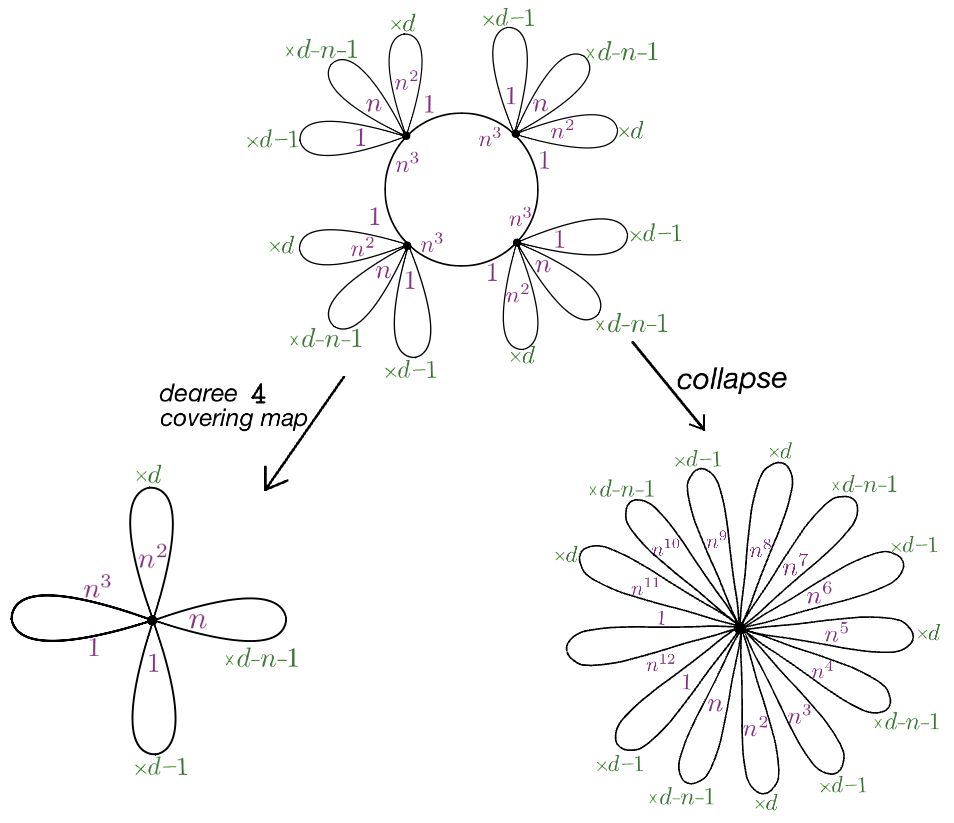}
    \caption{The labeled graph on the right represents an index $4$ subgroup $\Delta_{3,4}$ in $\Delta_3$.}
    \label{index l subgroup of Gamma}
\end{figure}

Since all entries of $\Vec{X}^{2l}(\Gamma_{1,l})$ are equal, $\Vec{X}^{2l}(\Gamma_{1,l})$ and any cyclic permutation of $\Vec{X}^{2l}(\Delta_{l,2})$ are linearly independent. Therefore, by Corollary \ref{dV(A)}, $\Gamma_{1,l}$ and $\Delta_{l,2}$ are not commensurable for all $l \geq 2$.

For $k, l\geq 3$, $\Vec{X}^{kl}(\Delta_{k,l})$ has $k-2$ consecutive equal terms, whereas $\Vec{X}^{lk}(\Delta_{l,k})$ has $l-2$ consecutive equal terms. Hence, $\Vec{X}^{kl}(\Delta_{k,l})$ and any cyclic permutation of $\Vec{X}^{kl}(\Delta_{l,k})$ are linearly independent, for $k, l\geq 3$. By the similar argument $\Vec{X}^{2k}(\Delta_{2,k})$ and any cyclic permutation of  $\Vec{X}^{2k}(\Delta_{k,2})$  are linearly independent for $k\geq 3$. This proves that for $k, l \geq 2$, $\Delta_{k,l}$ and $\Delta_{l,k}$ are commensurable if and only if $k=l$. Hence the same is true for $\Delta_k$ and $\Delta_l$.
\end{proof}

\subsection{Commensurability criterion for some GBS groups} \label{CRKZ_comm}
We conclude this section by giving a solution for the commensurability problem for a subclass of  GBS groups denoted by $\mathcal{C}_{n,l}$. This class of GBS group was introduced in \cite{MR4359918}.  For every $l \geq 1$, $n \geq 2$, $k \geq 2$, $0\leq a_1, a_2, \cdots, a_{k-1} \leq l-1$, denote by $A=A(n,l; a_1, a_2, \cdots, a_{k-1})$ the following labeled graph: it is bouquet of circles $e_1, \cdots, e_k$ with $\lambda(e_1)=1$, $\lambda(\overline{e}_1)=n^l$, and $\lambda(e_i)=\lambda(\overline{e}_i)=n^{a_{i-1}}$ for $2 \leq i \leq k$. 
The GBS group for this labeled graph is given by 
$$BS(1,n^l) \vee  BS(n^{a_1}, n^{a_1}) \vee  BS(n^{a_2}, n^{a_2}) \vee \cdots \vee BS(n^{a_{k-1}}, n^{a_{k-1}}).$$

 The solution to the isomorphism problem for groups in $\mathcal{C}_{n,l}$ is given in  of \cite[Theorem 5.1]{MR4359918}.

\begin{remark}
The vector associated in \cite{MR4359918} to a GBS group $G$ in $\mathcal{C}_{n,l}$ is $\frac{1}{2}\Vec{X}^l(G)$.
    
\end{remark}

\begin{theorem}\cite[Theorem 5.1]{MR4359918}
     Let $G_1$ and $G_2$ are two GBS groups in $\mathcal{C}_{n,l}$. Suppose $G_1$
     is the  GBS group defined by a labeled graph $A_1=A(n,l; a_1, a_2, \cdots, a_{k_1-1})$, and $G_2$ be the GBS group represented by a labeled graph $A_2=A(n,l; b_1, b_2, \cdots, b_{k_2-1})$. Then $G_1$ is isomorphic to $G_2$ if and only if 
     \begin{enumerate}
         \item $k_1=k_2$
    \item  $\Vec{X}^l(A_1)$ is a cyclic permutation of $\Vec{X}^l(A_2)$.
     \end{enumerate}
\end{theorem}

The following result provides a necessary and sufficient condition for two GBS groups in $\mathcal{C}_{n,l}$ to be commensurable.

\begin{figure}[h]
    \centering
    \includegraphics[width = 0.9 \textwidth]{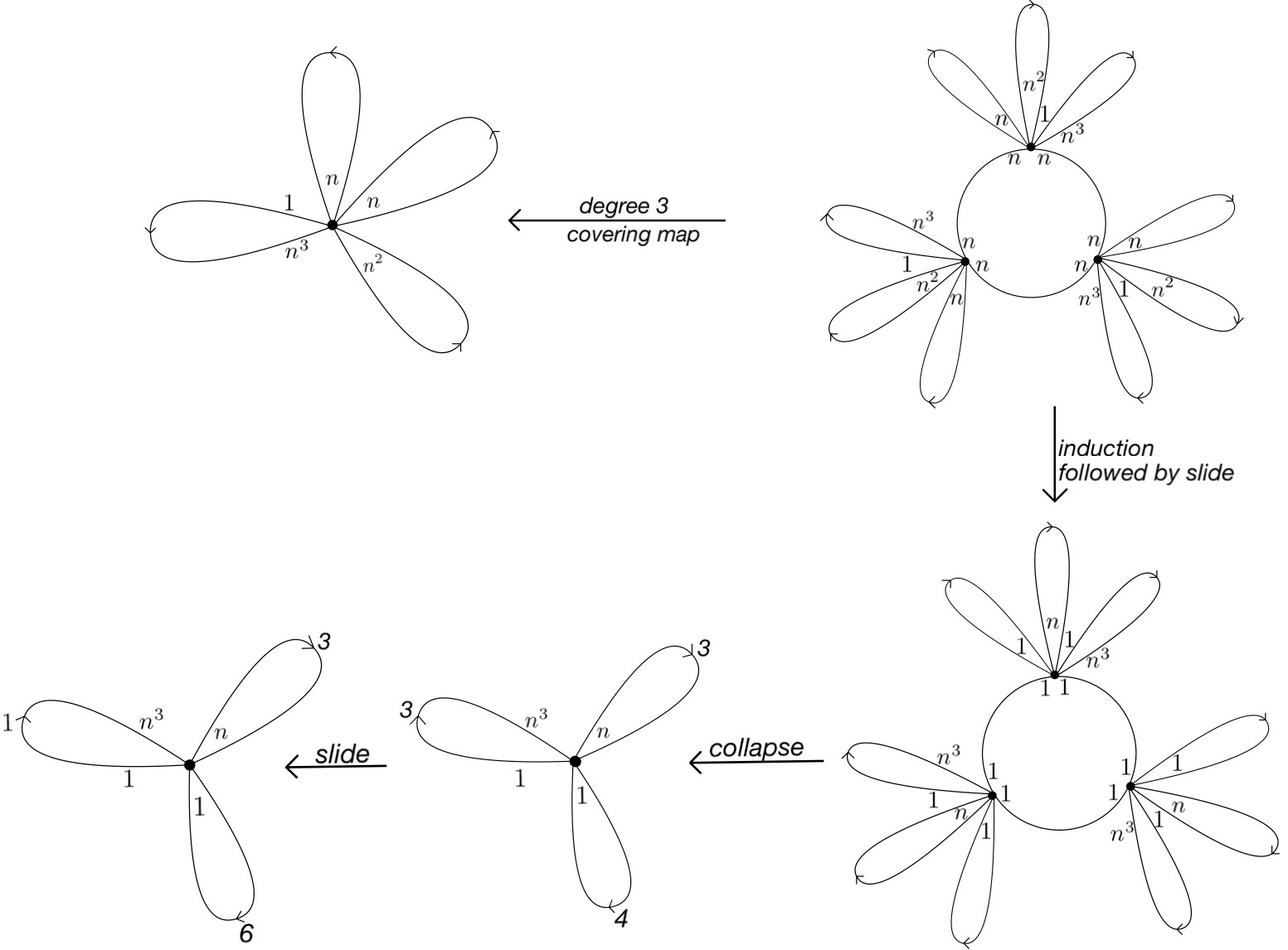}
    \caption{Example illustrating Theorem \ref{5, 15 thm} for $l=3$, $x^0=0$, $x^1=2$, $x^2=1$, $j_0=1$  and $c=3$. Here a number $z$ inside the petal means both ends of the petal have label $z$, and a number above a petal represents the multiplicity of that petal.}
    \label{fig:Bouquet of cicles eg}

\end{figure}

\begin{theorem}\label{5, 15 thm}
    Let $G_1$ and $G_2$ are two GBS groups in $\mathcal{C}_{n,l}$. Suppose $G_1$
     is the  GBS group defined by a labeled graph $A_1=A(n,l; a_1, a_2, \cdots, a_{k_1-1})$, and $G_2$ be the GBS group represented by a labeled graph $A_2=A(n,l; b_1, b_2, \cdots, b_{k_2-1})$.  Then $G_1$ is commensurable to $G_2$ if and only if $c_1\Vec{X}^l(A_1)$ is a cyclic permutation of $c_2\Vec{X}^l(A_2)$ for some $c_1, c_2 \in \mathbb{N} $.
\end{theorem}
\begin{proof}
\begin{figure}[h]
    \centering
    \includegraphics[width = 0.6 \textwidth]{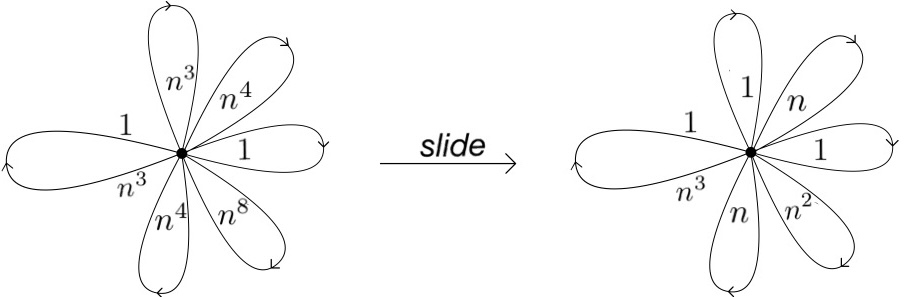}
    \caption{Slide moves on the bouquet of circles with $l=3$. Here $x^0=2$, $x^1=2$, and $x^2=1$.}
    \label{slide moves on A_i}

\end{figure}

The ``only if" direction follows from 
    Corollary \ref{dV(A)}. For the converse, by applying an induction move, we can assume that $c_1\Vec{X}^l(A_1)=c_2 \Vec{X}^l(A_2)$. Perform slide moves on labeled graphs $A_i$ for $i=1,2$ (see Figure\ref{slide moves on A_i} for an example) to obtain the bouquet of circles whose GBS group which is isomorphic to $G_i$ is the following:
     $$ BS(1,n^l) \vee \bigvee_{x^0_i} BS(1, 1) \vee \bigvee_{x^1_i} BS(n,n) \vee \cdots \vee \bigvee_{x^{l-1}_i} BS(n^{l-1}, n^{l-1})$$

     Therefore, the graph $A_i$ is in the same deformation space as the graph  $A_i'$, which is defined by a bouquet of circles with $x^j_i$ circles each with label $n^j$ on both ends for each $j$, and one circle with label $1$ on one end and $n^l$ on the other end.  One can see that $\Vec{X}^l(A_i)=\Vec{X}^l(G_i)=2(x_i^0, 
    x_i^1, \cdots, x_i^{l-1})$. 
   
 Suppose $j_0 \leq l-1$ is the smallest number such that $x^{j_0}_i \neq 0$, and   $x^{j}_i = 0$ for all $0 \leq j < j_0$. Let $B_i'$ denote a $c_i$ sheeted topological covering of $A_i'$ that unwinds a loop in $A_i'$ with labels $n^{j_0}$ on both ends, into a cycle of length $c_i$ (Figure  \ref{fig:Bouquet of cicles eg} illustrate an example of this). We can apply an induction move to each vertex of $B_i'$ to make the labels on the cycle $n^l$. By performing slide moves to each vertex of the resulting labeled graph, we can get all labels on the cycle to be $1$.
Now, by applying $c_i-1$ collapse moves, we obtain a bouquet of circles. We can adjust the petal labels by applying slide moves, resulting in $c_1x^j_1=c_2x^j_2$ petals with labels $n^{j-j_0+1}$ on both ends for $ 0 \leq j \leq l-1$, along with one circle having a label $1$ on one end and $n^l$ on the other end. 
\\
Thus, the labeled graphs $B_1'$ and $B_2'$ are related by slide moves to the same labeled graph. Consequently, they represent isomorphic finite index subgroups of $G_1$ and $G_2$, respectively. This completes the proof of the theorem.   
    
\end{proof}

\section{Reduced graphs with no proper plateau}\label{eg using plateau}
In the following section, we will give examples of incommensurable lattices in Aut$(X_{dm,dn})$, when $\gcd(m,d)=\gcd(n,d)=1$ for $ m,n > 1$, and
there is a prime $q < d$ which divides either $m$ or $n$. Without loss of generality, we can assume $q$ divides $m$ (hence $\gcd(m,q)=q$ and $\gcd(d,q)=1$).

Let $D$ be the largest number dividing $d-q$ which is coprime to both $m$ and $n$. Let $d-q=ID$.

\begin{figure}[h]
    \centering
    \includegraphics[width=0.45\linewidth]{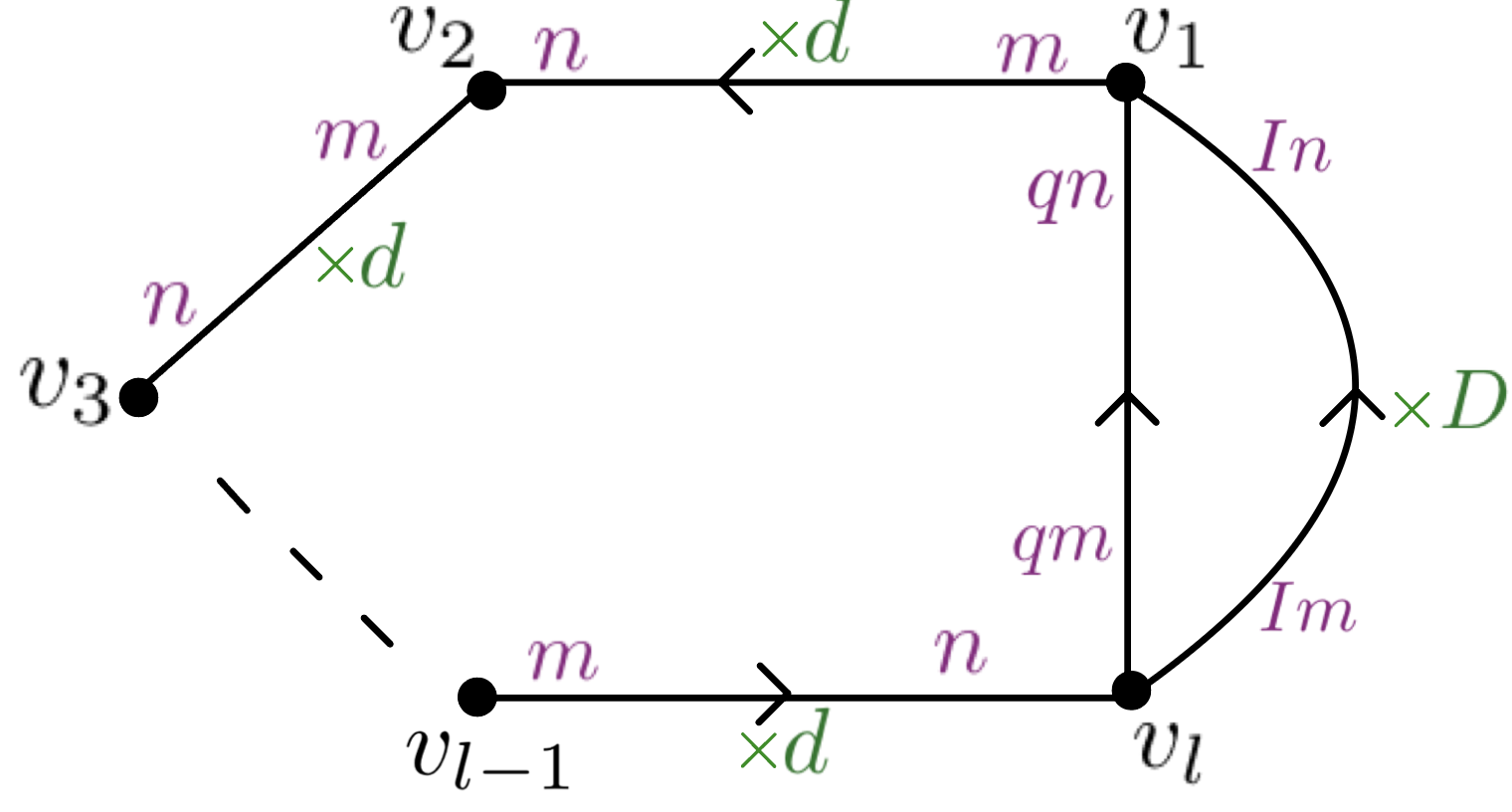}
    \caption{Graph $B_l$ for $l \geq 2$, defining lattices in Aut$(X_{dm,dn})$ when $gcd(m,d) = gcd(n,d) = 1$ and a prime $q < d$ divides $m$.}
    \label{depth profile gammas}
\end{figure}

\textbf{The lattice {$\mathbf{\Lambda_l}$} for $\mathbf{l \geq 2}$.} Consider the lattice $\Lambda_l$ defined by the directed labeled graph $(L_l, \lambda_l)$ given in Figure (\ref{depth profile gammas}). It is a graph with $l$ vertices $v_1, v_2, \cdots , v_l$ and $d(l-1)+D+1$ directed edges. There are $d$ directed edges from $v_{i}$ to $v_{i+1}$ for $1 \leq i \leq k-1$ with initial label $m$ and terminal label $n$. One directed edge from $v_l$ to $v_1$ with initial label $qm$, and terminal label $qn$. Lastly, there are $D$ edges from $v_l$ to $ v_1$ with initial label $Im$, and terminal label $In$.

Recall the graph $(B_1, \mu_1)$ defined by Figure \ref{fig:B_1} and its fundamental group $\Gamma_1$. Since $m, n > 1$, the labeled graphs $(B_1, \mu_1)$ and $(L_l, \lambda_l)$ are reduced for all $l \geq 1$. The next two lemmas imply that these labeled graphs do not contain any proper $p$--plateau. Hence, the finite index subgroups of $\Gamma_1$  and $\Lambda_l$ are given by topological coverings of $B_1$ and $L_l$ for $l \geq 2$. By analyzing these covers, we will show that any two groups in $\{\Gamma_1$, $\Lambda_l: l \geq 2\}$  are abstractly incommensurable.

\begin{remark}\label{no_plateau} For the labeled graph $(A, \lambda)$ the following holds which can also be found in \cite{levitt1}:
    \begin{enumerate}
    
        \item For a vertex $v\in V(A)$, $\{v\}$ is a  $p$--plateau if and only if $p$ divides every label at $v$. 
        \item 
        Let $P\subseteq A$ be a $p$--plateau and $e \in E(A)$. Then the following holds:
        \begin{enumerate}
            \item If $e \in E(P)$, then $p \nmid \lambda(e), \lambda(\overline{e})$.
           
            \item If $\partial_0(e), \partial_1(e) \in V(P)$ but $e \not \in E(P)$, then $p \mid \lambda(e), \lambda(\overline{e})$.
            
            \item If $\partial_0(e) \in P$ and $\partial_1(e) \notin V(P)$, then $p \mid \lambda(e)$.
           
            \item If $\partial_0(e), \partial_1(e) \notin P$, then there is no restriction.
        \end{enumerate}
        \item If $P\subseteq A$ is a $p$--plateau for prime $p$ not dividing any label of $A$, then $P=A$.
    \end{enumerate}
\end{remark}

\begin{lemma}
    $(B_1, \mu_1)$ does not contain a proper plateau. 
\end{lemma}
\begin{proof}
    $\{u\}, \{v\} \subset B_1$ are not $p$--plateaus for any prime number $p$ by Remark \ref{no_plateau}(1) and the fact that  $\gcd(m,n)=1$. 

If $e\in E(P)$ for some $p$--plateau $P\subseteq B_1$, then by Remark \ref{no_plateau}(2a), $p\nmid \mu_1(e), \mu_1(\overline{e})$. Therefore $p$  does not divide $m$ or $n$ and by Remark, \ref{no_plateau}(3) $P=B_1$.
\end{proof}

\begin{lemma}
    $(L_l, \lambda_l)$ does not contain a proper plateau for $l\geq 2$. 
\end{lemma}
\begin{proof}
Let $P \subseteq L_l$ be a $p$--plateau for a prime number $p$. If $p$ is coprime to all $m,n,q,I$,  then by Remark \ref{no_plateau}(3), $P=L_l$, hence $P$ is not a proper  $p$--plateau. We will show that if $p$ divides any of the numbers $m,n,q, I$,  then $P$ is a null graph. 

 If $p \mid q$, then $p=q$ and $q \mid m$. Since $\gcd(m,n)=1$ and $q \mid m$, we have $q\nmid n$. We claim that $q\nmid I$.
Assuming $q\mid I$, we will have the following sequence of implications, contradicting the fact that $\gcd(d,q)=1$;
$$q \mid I \implies q \mid DI \implies q \mid d-q \implies q \mid d.$$

$P \neq \{v_i\}$ for $1 \leq i \leq l$ by Remark \ref{no_plateau}(1) and the fact that $q\nmid I, n$.
By Remark \ref{no_plateau} (2a), $P$ doesn't contain any edge since $q \mid m$. Therefore $P$ is a null graph.

 If $p \mid m$, then $p\nmid n$ as $\gcd(m,n)=1$. If $p=q$ then we have already shown that $L_l$ does not contain proper $q$-plateau. If $p \neq q$, then $P\neq \{v_i\}$ by Remark \ref{no_plateau}(1) and the fact that $p\nmid n, q$. 
Moreover, the edges are not contained in $P$ as $p \mid m$.

  If $p\mid n$, then $p \nmid m$ and $p \nmid q$.
$P\neq \{u\}, \{v\}$ due to  Remark \ref{no_plateau}(1) and the fact that $p\nmid m, q$. 
The edges are not contained in $P$ as $p \mid n$.

We have seen that $L_l$ does not contain a proper $p$--plateau for $p$ dividing $m,n$ or $q$. Therefore, assume $p\nmid m,n,q$. Recall that $D$ was chosen to be the largest number such that $D \mid (d-q)$ and $\gcd(D,n)=\gcd(D,m)=1$. If $p \mid I$, then $d-q=DI$ implies $pD \mid (d-q)$. Also,  $\gcd(D,n)=\gcd(D,m)=1=\gcd(p,m)=\gcd(n,p)$  implies $\gcd(pD,n)=\gcd(pD,m)=1$, contradicting the choice of $D$. 
\end{proof}

\begin{proposition}\label{mainprop_case2}
    For $i\in \{1,2\}$, let $G_i$ be the GBS groups represented by a labeled graphs $(A_i,\lambda_i)$ with no proper plateau and such that $\lambda(e_i) \neq \pm1$ for any $e_i \in E(A_i)$. Furthermore, suppose $q(G_i) \cap \mathbb{Z} = \{1\}$. Then if $G_1$ and $G_2$ are commensurable, then  $\frac{|V(A_1)|}{|V(A_2)|}=\frac{|E(A_1)|}{|E(A_2)|}$.

\end{proposition}
\begin{proof}
    Since  $G_1$ and $G_2$ are commensurable, there exist admissible covers $(\Tilde{A_i}, \Tilde{\lambda_i})$ of $(A_i, \lambda_i)$ representing isomorphic GBS groups. 
    
    Since $(A_i, \lambda_i)$ does not contain a proper plateau, by Proposition \ref{Levitt}, $\Tilde{A_i}$ is a topological cover of $A_i$ of some degree $d_i$. Therefore, 
\begin{equation}\label{equa2}
    |V(\Tilde{A_i})|=d_i |V(A_i)| \text{ and }  |E(\Tilde{A_i})|=d_i |E(A_i)|.
\end{equation}

Since $\lambda(e_i) \neq \pm 1$ for any $e_i \in E(A_i)$, and $\Tilde{A_i}$ is  a topological cover of $A_i$, we have $\lambda(\Tilde{e}_i) \neq \pm 1$ for any $\Tilde{e}_i \in E(\Tilde{A}_i)$. Therefore, $\Tilde{A_i}$ is a reduced graph. 
By Proposition \ref{reduced graph related by slide move}, $(\Tilde{A_1}, \Tilde{\lambda_1})$ and $(\Tilde{A_2}, \Tilde{\lambda_2})$ represent isomorphic GBS groups if and only if they are related by slide moves. Since slide moves do not change the numbers of vertices and edges in a graph,  we have
\begin{equation}\label{equa1}
|V(\Tilde{A_1})|=|V(\Tilde{A_2})| \text{ and } |E(\Tilde{A_1})|=|E(\Tilde{A_2})|.
\end{equation}

    From equations (\ref{equa1}) and (\ref{equa2}),  we get the following equality;
   $$\frac{|V(A_1)|}{|V(A_2)|}=\frac{d_2}{d_1}=\frac{|E(A_1)|}{|E(A_2)|}.$$
\end{proof}

\begin{corollary}\label{d neq D+1}
    The GBS groups $\Gamma_1$ and $\Lambda_l$ are not abstractly commensurable for $l\geq 2$. 
\end{corollary}
\begin{proof}

Note that $B_1$ and $L_l$ is reduced graph for all $l \geq 2$ as $m,n > 1$. Also, since $\gcd(m,n)=1$, it follows that,  $q(\Gamma_k)\cap \mathbb{Z}=\{1\}$. 

Now, assume $\Gamma_1$ and $\Gamma_l$ are commensurable groups for $l \geq 2$. By Proposition \ref{mainprop_case2} we have
    $$\frac{|V(B_1)|}{|V(L_l)|}=\frac{|E(B_1)|}{|E(L_l)|}$$
    which is equivalent to 
     $$\frac{2}{l}=\frac{2d}{d(l-1)+D+1}$$
   Rearranging and simplifying this equation yields $d=D+1$. Since $(L_l, \lambda_l)$ is a uniform lattice in Aut$(X_{dm,dn})$, by Proposition \ref{prop4.6} and the fact that $I \geq 1$, we get the contradiction
    \begin{align*}
    dm&=\sum_{e \in E^+_0(v_l)} \mu_l(e)\\
    &=qm+mDI\\
    &\geq qm+mD \\
    &= qm+(d-1)m\\
&\geq 2m+(d-1)m\\
&= (d+1)m
\end{align*}
where the last inequality follows from the fact that $q$ is a prime number, and hence $q\geq 2$.
\end{proof}

\begin{corollary}
    The GBS groups $\Lambda_k$ and $\Lambda_l$ are not abstractly commensurable for   $k,l \geq 2$ and $k \neq l$ . 
\end{corollary}
\begin{proof}
 Assume  $\Lambda_k$ and $\Lambda_l$ are commensurable groups. Then, by Proposition \ref{mainprop_case2}, we have the following statements:
    $$\frac{|V(L_l)|}{|V(L_k)|}=\frac{|E(L_l)|}{|E(L_k)|}$$
    
    $$\frac{l}{k}=\frac{d(l-1)+D+1}{d(k-1)+D+1}$$
      $$l(k-1)d+l(D+1))=k(l-1)d+k(D+1))$$
       $$l(D-d+1)=k(D-d+1).$$
       Since $d \neq D+1$ (by the same argument as in the proof of \ref{d neq D+1}), it follows that $k=l$.  This completes the proof of the corollary.
      
\end{proof}

Finally, we conclude this paper by giving necessary and sufficient conditions for the cell complex $X_{dm,dn}$ to satisfy Leighton's Property.

\begin{theorem}
    The Baumslag-Solitar complex $X_{dm,dn}$ has the Leighton property if and only if $m$ and $n$ have no divisor less than or equal to $d$.
\end{theorem}
\begin{proof}
When  $m$ or $n$ has a divisor less than or equal to $d$, sections \ref{eg using DP}, \ref{eg using crkz}, and \ref{eg using plateau} provide examples of incommensurable lattices in Aut$(X_{dm,dn})$. This proves the forward direction of the Theorem. 
The converse direction is proved in Theorem \ref{one direction of 1.4}
\end{proof}

\bibliographystyle{alpha}

\bibliography{Paper_24}

\begin{thebibliography}{CRKZ21}

\bibitem[Bas93]{bass}
Hyman Bass.
\newblock Covering theory for graphs of groups.
\newblock {\em J. Pure Appl. Algebra}, 89(1-2):3--47, 1993.

\bibitem[BL01]{basslubotzky}
Hyman Bass and Alexander Lubotzky.
\newblock {\em Tree lattices}, volume 176 of {\em Progress in Mathematics}.
\newblock Birkh\"{a}user Boston, Inc., Boston, MA, 2001.
\newblock With appendices by Bass, L. Carbone, Lubotzky, G. Rosenberg and J. Tits.

\bibitem[BM00]{MR1839489}
Marc Burger and Shahar Mozes.
\newblock Lattices in product of trees.
\newblock {\em Inst. Hautes \'Etudes Sci. Publ. Math.}, (92):151--194, 2000.

\bibitem[BS22]{MR4464467}
Martin~R. Bridson and Sam Shepherd.
\newblock Leighton's theorem : extensions, limitations and quasitrees.
\newblock {\em Algebr. Geom. Topol.}, 22(2):881--917, 2022.

\bibitem[CRKZ21]{MR4359918}
Montserrat Casals-Ruiz, Ilya Kazachkov, and Alexander Zakharov.
\newblock Commensurability of {B}aumslag-{S}olitar groups.
\newblock {\em Indiana Univ. Math. J.}, 70(6):2527--2555, 2021.

\bibitem[DK23]{MR4545211}
Natalia~S. Dergacheva and Anton~A. Klyachko.
\newblock Small non-{L}eighton two-complexes.
\newblock {\em Math. Proc. Cambridge Philos. Soc.}, 174(2):385--391, 2023.

\bibitem[For06]{MR2276234}
Max Forester.
\newblock Splittings of generalized {B}aumslag-{S}olitar groups.
\newblock {\em Geom. Dedicata}, 121:43--59, 2006.

\bibitem[For24]{forester2022incommensurable}
Max Forester.
\newblock Incommensurable lattices in {B}aumslag-{S}olitar complexes.
\newblock {\em J. Lond. Math. Soc. (2)}, 109(3):Paper No. e12879, 31, 2024.

\bibitem[Lei82]{leighton}
Frank~Thomson Leighton.
\newblock Finite common coverings of graphs.
\newblock {\em J. Combin. Theory Ser. B}, 33(3):231--238, 1982.

\bibitem[Lev07]{levitt-auto}
Gilbert Levitt.
\newblock On the automorphism group of generalized {B}aumslag-{S}olitar groups.
\newblock {\em Geom. Topol.}, 11:473--515, 2007.

\bibitem[Lev15]{levitt1}
Gilbert Levitt.
\newblock Generalized {B}aumslag-{S}olitar groups: rank and finite index subgroups.
\newblock {\em Ann. Inst. Fourier (Grenoble)}, 65(2):725--762, 2015.

\bibitem[Neu10]{MR2727669}
Walter~D. Neumann.
\newblock On {L}eighton's graph covering theorem.
\newblock {\em Groups Geom. Dyn.}, 4(4):863--872, 2010.

\bibitem[She22]{MR4506536}
Sam Shepherd.
\newblock Two generalisations of {L}eighton's theorem.
\newblock {\em Groups Geom. Dyn.}, 16(3):743--778, 2022.
\newblock With an appendix by Giles Gardam and Daniel J. Woodhouse.

\bibitem[SW24]{MR4748179}
Emily~R. Stark and Daniel~J. Woodhouse.
\newblock Action rigidity for free products of hyperbolic manifold groups.
\newblock {\em Ann. Inst. Fourier (Grenoble)}, 74(2):503--544, 2024.

\bibitem[Wis96]{MR2694733}
Daniel~T. Wise.
\newblock {\em Non-positively curved squared complexes: {A}periodic tilings and non-residually finite groups}.
\newblock ProQuest LLC, Ann Arbor, MI, 1996.
\newblock Thesis (Ph.D.)--Princeton University.

\bibitem[Wis07]{MR2341837}
Daniel~T. Wise.
\newblock Complete square complexes.
\newblock {\em Comment. Math. Helv.}, 82(4):683--724, 2007.

\bibitem[Woo21]{MR4243770}
Daniel~J. Woodhouse.
\newblock Revisiting {L}eighton's theorem with the {H}aar measure.
\newblock {\em Math. Proc. Cambridge Philos. Soc.}, 170(3):615--623, 2021.

\bibitem[Woo23]{MR4647679}
Daniel~J. Woodhouse.
\newblock Leighton's theorem and regular cube complexes.
\newblock {\em Algebr. Geom. Topol.}, 23(7):3395--3415, 2023.

\end{thebibliography}
\end{document}